\documentclass[amstex,amstext,amsfonts,epsf,12pt,color,a4paper]{article}%
\usepackage{epsfig,latexsym,amsfonts,amsmath,amssymb,dsfont}
\usepackage{graphicx}
\usepackage{bm}
\usepackage{color}

\newtheorem{theorem}{Theorem}[section]

\newtheorem{lemma}[theorem]{Lemma}
\newtheorem{proposition}[theorem]{Proposition}

\newtheorem{summary}[theorem]{Summary}

\def\be{\begin{eqnarray*}}
\def\ee{\end{eqnarray*}}
\def\ben{\begin{eqnarray}}
\def\een{\end{eqnarray}}
\def\IntO{\int\limits_\Omega}
\def\mean#1{\left\{#1\right\}}
\def\wh{\widehat}
\def\v{{\bf v}}
\def\f{{\bf f}}

\def\wt{\widetilde}
\def\cM{{\mathcal M}}
\def\cT{{\mathcal T}}
\def\be{\begin{eqnarray*}}
\def\ee{\end{eqnarray*}}
\def\ben{\begin{eqnarray}}
\def\een{\end{eqnarray}}
\def\IntO{\int\limits_\Omega}
\def\mean#1{\left\{\left\vert #1\right\vert\right\}}
\def\wh{\widehat}
\def\v{{\bf v}}
\def\f{{\bf f}}
\def\abig{{\rm a}}

\makeatletter\@addtoreset{equation}{section}\makeatother
\makeatletter\@addtoreset{figure}{section}\makeatother
\makeatletter\@addtoreset{table}{section}\makeatother

\newcommand{\Order}[1]{\mathcal{O}\left(#1\right)}

\DeclareMathOperator*{\diag}{diag}

\makeatletter\@addtoreset{equation}{section}\makeatother
\makeatletter\@addtoreset{figure}{section}\makeatother
\makeatletter\@addtoreset{table}{section}\makeatother
\def\mean#1{\left\{\,#1\,\right\}}

\def\Rd{{\mathbb R}^d}

\begin{document}
\title
{
A fast iteration method for solving elliptic problems with
quasiperiodic coefficients 
}

\author{Boris N. Khoromskij\thanks{Max Planck Institute for
        Mathematics in the Sciences, Inselstr.~22-26, D-04103 Leipzig,
        Germany ({\tt bokh@mis.mpg.de})}
        \and Sergey I. Repin\thanks{V.A. Steklov Institute of Mathematics, 191011, 
        Fontanka 27, St.Petersburg, Russia ({\tt repin@pdmi.ras.ru})}}


\date{}

\maketitle
\begin{flushright}\it Dedicated to 70th jubilee of Prof. Yu. A. Kuznetsov
\end{flushright}
\vskip8pt

\begin{abstract}
{}\small
The paper suggests a preconditioning type method for fast solving of  elliptic equations 
with  oscillating  quasiperiodic coefficients $A_\epsilon$ specified by  the
small parameter $\epsilon>0$.
We use an iteration method generated by an elliptic operator,
associated with a certain simplified (e.g., homogenized) problem. On each step of this
procedure it is required to solve an auxiliary elliptic boundary value problem
with non--oscillating coefficients $A_0$.
All the information related to complicated coefficients of the original differential problem
is encompasses in the linear functional, which forms the right hand side of the  auxiliary problem.   
Therefore, explicit inversion of the original operator associated with oscillating
coefficients is avoided. The only operation used instead
is multiplication of the operator by a vector (vector function), which 
can be efficiently performed due to the low-rank QTT tensor 
operations with the rank parameter controlled 
by the given precision $\delta >0$ independent on the parameter $\epsilon$. 
In the first part of the paper, we establish sufficient conditions that guarantee
convergence of  the iteration method and deduce
explicit estimates of the contraction factor, which
are expressed in terms of $A_\epsilon$ and $A_0$.
Moreover, we deduce two--sided a posteriori error estimates that  do not use 
$A^{-1}_\epsilon$ and provide guaranteed two sided bounds of the distance to the exact solution
of the original problem for any step of the iteration process. The second part
is concerned with  realisations of the iteration method.
For a wide class of oscillating coefficients, we obtain  sharp QTT 
rank estimates for the stiffness matrix in tensor representation.
In practice, this leads to the logarithmic complexity scaling of the approximation and solution
process in both the FEM grid-size, and  $O(\vert\log\epsilon\vert)$ cost in terms of $\epsilon$.
Numerical tests in 1D confirm the logarithmic complexity scaling of our method applied
to a class of  complicated quasiperiodic coefficients.
\end{abstract}



\section{Introduction}
\label{sec:Intr_problem_set}

Partial differential equations with  oscillating coefficients
often arise in various models in natural sciences, including quantum chemistry and material
sciences, as well as in engineering applications.
Numerical analysis of problems with periodical coefficients is often performed by
geometric homogenization methods, which provide efficient approximations of structures
with very large amount of cells of periodicity
(see, e.g., \cite{Bakhvalov,Bensoussan,Cioranescu,Jikov,GlorOtto:12}). 
We consider a wider class of problems where either the amount of cells is significant
but not large enough to ignore modeling errors generated
by homogenized models or periodicity has a more complicated form. Numerical analysis of such problems is
faced with 
 several challenging problems. The main three of them
 are as follows: (a) creation of a robust numerical 
 method able to
 construct a sequence converging to the exact solution
 by means of using finite element approximations on regular (quasiregular) meshes; (b) guaranteed a posteriori
 estimates of the distance between the exact solution
 of a boundary value problem with highly oscillating
 coefficients and an approximation; (c) construction
 an efficient solver based on suitable
preconditioning of the respective discrete system.

\begin{figure}[htbp]
 \includegraphics[width=5.5cm]{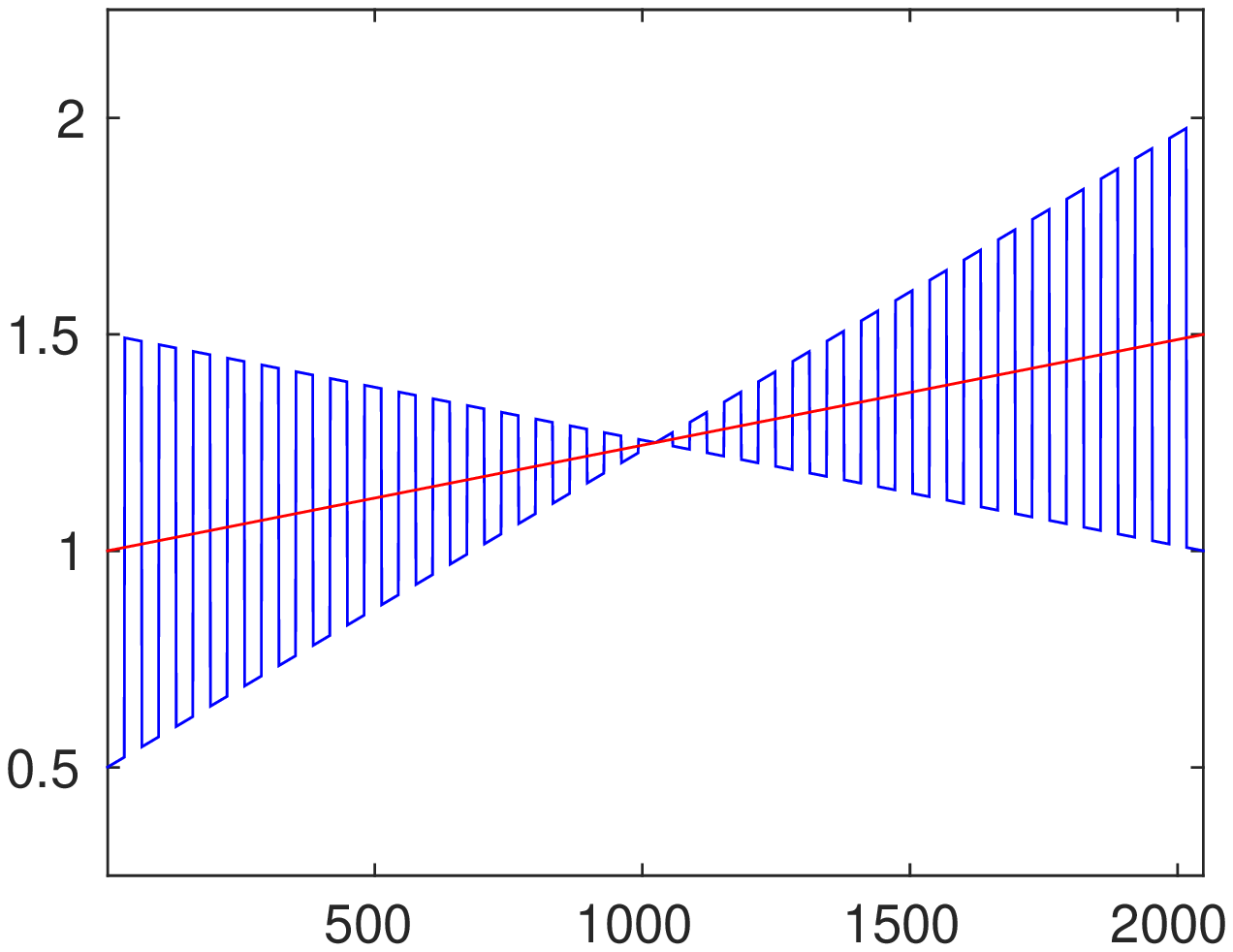}
\includegraphics[width=5.5cm]{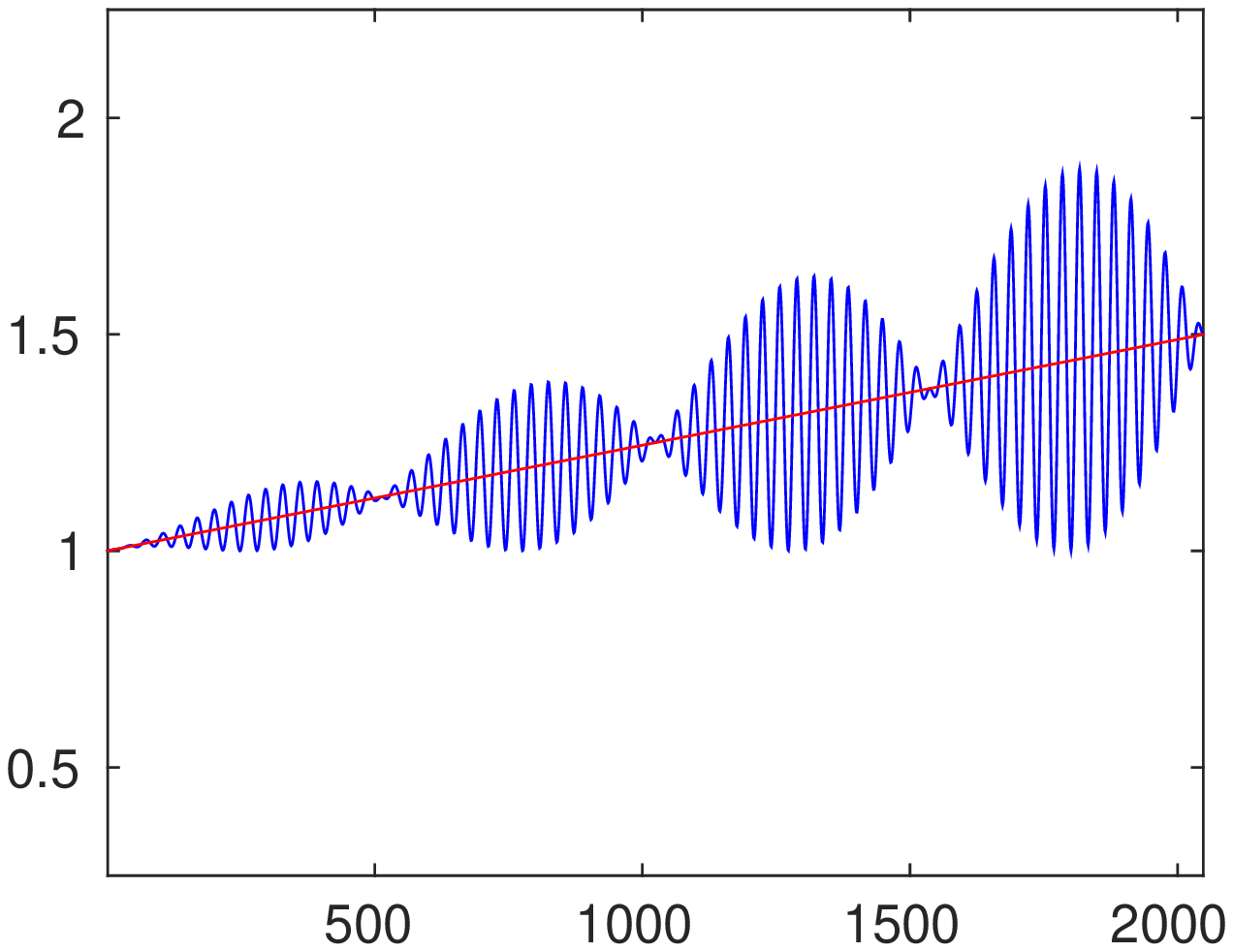}
\includegraphics[width=5.3cm]{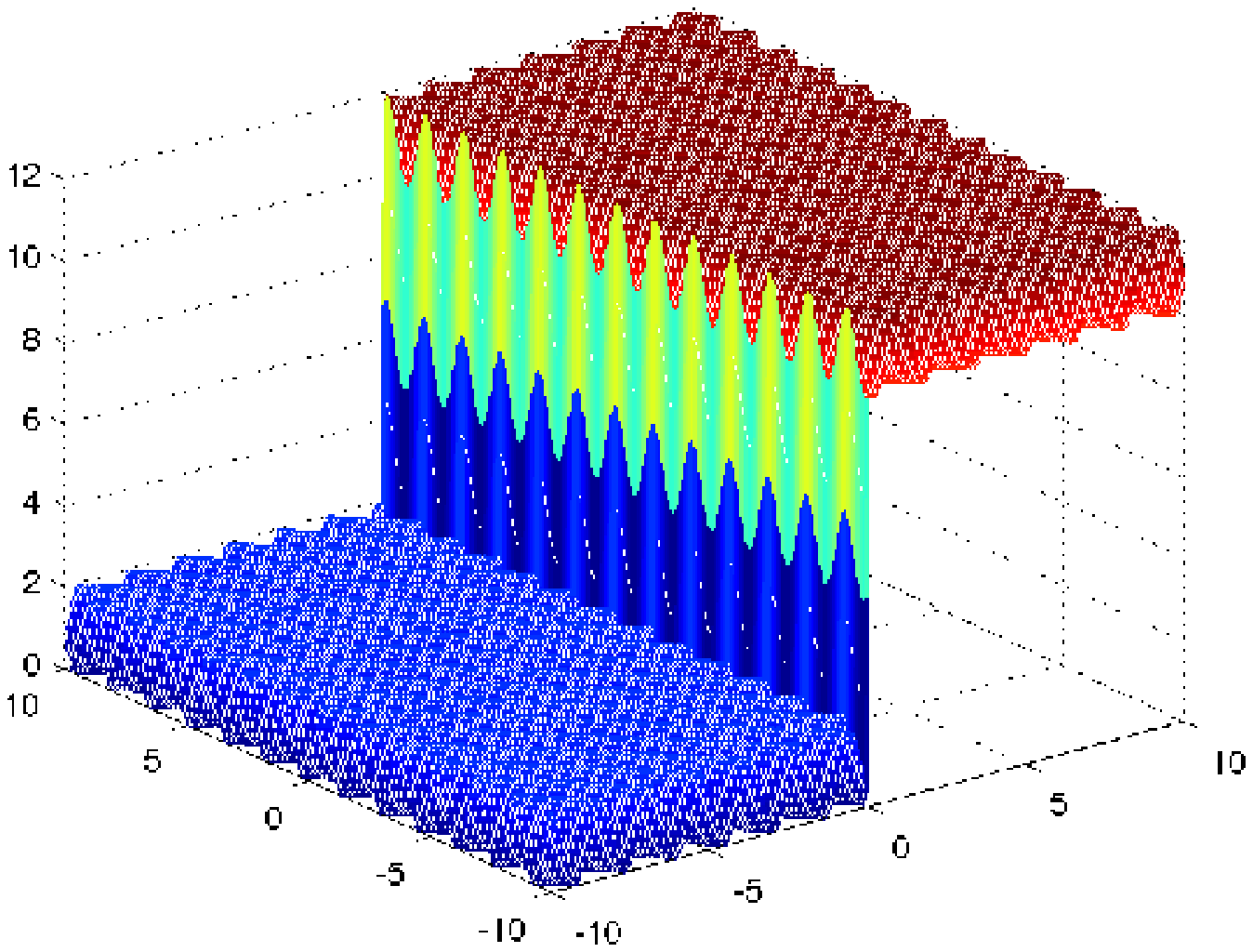}
\caption{ Examples of modulated periodic and piecewise periodic coefficients in 1D.
}
\label{fig:3DPeriodStruct2}
\end{figure}

In this paper, we suggest an approach that solves (a)--(c) 
for a class of elliptic problems with quasi-periodic coefficients. We discuss the basic ideas
with the paradigm of the model second-order elliptic problem,
but it is clear that they can be extended to many
other elliptic and parabolic type equations with quai--periodic coefficients.
Consider the problem
\begin{equation}
\label{1.1}
-\operatorname{div}\left(  A_{\epsilon}(x)\nabla u_{\epsilon}\right)
=f\quad\text{in }\Omega,
\qquad u_\epsilon=0\quad{\rm on}\,\partial\Omega,
\end{equation}
where $f\in L^{2}\left(  \Omega\right)$, $\Omega=(0,1)^d$ ($d=1,2,3$),
with homogeneous Dirichlet boundary conditions,  where a small parameter 
$\epsilon >0$ is a small parameter characterizing oscillations,
and  $A_{\epsilon}(x)$ is a matrix with quasiperiodic coefficients.
  We assume that
 ${A_\epsilon}\in L^{\infty
}(\Omega,\,\mathbb{M}_{\mathrm{sym}}^{d\times d})$ (here and later on $\,\mathbb{M}_{\mathrm{sym}}^{d\times d}$ denotes the set 
of symmetric $\,d\times d-$ matrices) and
 \begin{equation}
\label{1.2}
\lambda^\epsilon_\ominus |\zeta|^2\,\leq A_\epsilon(x)\zeta\cdot\zeta\,\leq \lambda^\epsilon_\oplus|\zeta|^2,
\qquad\forall\,\zeta\in{\mathbb R}^d,\;x\in \Omega,
\end{equation}
where $\lambda^\epsilon_\ominus$ is a positive constant,
so that the problem is well posed and the corresponding generalized
solution $\,u_{\epsilon}\in H_{0}^{1}(\Omega)\,$ is defined by the relation
\begin{equation}
{\abig}_\epsilon(u_\epsilon,w)
=\,(f,w)_\Omega
\in H_{0}^{1}\left(  \Omega\right),   
\label{1.3}%
\end{equation}
where
 $$
 \abig_\epsilon(u,w):=\int_{\Omega}A_{\epsilon}\nabla u_\epsilon\cdot\nabla w dx\quad{\rm  and}\quad
  (f,w)_\Omega:=\int_{\Omega}f w dx.
  $$
Entries of  $A_\epsilon$  may depend on $x$ in a very complicated
 way, see some examples depicted in Fig. \ref{fig:3DPeriodStruct2}.
Therefore,  the problem (\ref{1.1})
 may be very  difficult from the viewpoint of quantitative analysis. The level of complexity
 can be roughly estimated by the parameters $\kappa:=\frac{\lambda^\epsilon_\ominus}{\lambda^\epsilon_\oplus}$ and  $\epsilon$. 
 If both of them are  very small,
 then serious difficulties will arise in  approximation  methods
 and in numerical
solution of the corresponding linear systems
(which may have very large dimensions and huge condition
numbers). For such type problems, getting guaranteed and efficient a posteriori error 
estimates may be a highly difficult problem as well.

Our goal is to justify a numerical method for computing
successful approximations of $u_\epsilon$ which is based
on solving a simpler problem associated with the bilinear
form $\abig_0(u,w)=\IntO A_0\nabla u\cdot\nabla w\,dx$.
Coefficients of the matrics $A_0$ are much more regular than coefficients of $A_\epsilon$ and do not have rapid
oscillations. It is assumed that $A_0$ satisfies the condition
\ben
\label{1.5}
\lambda^0_\ominus|\zeta|^2\quad\leq &A_0\zeta\cdot\zeta& \leq \quad\lambda^0_\oplus
|\zeta|^2\qquad
\forall\zeta\in{\mathbb R}^d,\;x\in \Omega
\een
with  positive constants $\lambda^0_\ominus$ and $\lambda^0_\oplus$. Then,
there exist positive constants $\lambda_1$
and $\lambda_2$ such that
\ben
\label {1.10}
\lambda_1 A_0\zeta\cdot\zeta
\leq \,
A_\epsilon\zeta\cdot\zeta\,\leq \lambda_2 A_0\zeta\cdot\zeta\qquad\forall\zeta\in \Rd,\quad
x\in \Omega.
\een

The  homogenization theory
 suggests a suitable form of $\abig_0(u,w)$ for perfectly
periodical structures, where 
$\Omega$ is a collection of self-similar cells 
$\Pi^\epsilon_i$, $i=1,2,...,L$ 
and  the cell size $\epsilon$
is very small (in comparison with the ${\rm diam}\Omega$). 
In this case, for any ${x}\in\Pi_{{i}}^{\epsilon}$ the matrix is defined
by the relation
$
A_{\epsilon}({x}):=\widehat{A}({ y })\in L^{\infty
}(\widehat{\Pi},\,\mathbb{M}_{\mathrm{sym}}^{d\times d})
$, where
${y}=
\frac{{x-\zeta_{i}}%
}{\epsilon}$, $\zeta_i$ is the "cell centre", and $y$ is the  cartesian coordinate system
associated with the "reference" cell
$\,\widehat{\Pi}$. 
An approximation of $\,u_{\epsilon}\,$ is constructed by a special procedure.
First, for $\,k=1,2,...,d$
we find the solutions $\,N_k\,$ of ``cell problems''
\begin{equation}
\label{1.6}%
 {\rm div}(\,\widehat{A}\,\nabla\,N_{k})=({\rm div}%
\,\,\widehat{A})_{k}\quad\mathrm{in}\quad\widehat{\Pi},
\end{equation}
which
satisfy the the periodic boundary conditions and the mean value condition
\begin{equation*}
\mean{N_k}_{\widehat \Pi}:=
\frac{1}{|\widehat \Pi|}\int_{\widehat{\Pi}}N_{k}=0.
\label{c_problem}%
\end{equation*}
Then, we define the  matrix
$
A_{0}=\;\mean{\widehat{A}\,(\,I-\nabla\boldsymbol{N})}_{\wh \Pi},
$
where
$\boldsymbol{N}=\{N_1,N_2,...,N_d\}$. 
The homogenized problem is to find $\,u_{0}\in
H_{0}^{1}(\Omega)\,$ such that
\begin{equation}
\abig_0(u_0,w)\,=(f,w)_\Omega\qquad\forall
w\in H_{0}^{1}(\Omega),
 \label{1.7}%
\end{equation}
where $\abig_0(u_0,w):=\int_{\Omega}\,A_{0}\nabla u_{0}\cdot\nabla w dx$.
This problem is much simpler than the original one. The function $u_0$ approximates
$\,u_{\epsilon}$ in a weak sense (see, e.g., \cite{Bensoussan}),
\begin{equation*}
u_{\epsilon}\rightarrow u_{0}
\quad\mathrm{in}\;L^{2}(\Omega)\qquad
\mathrm{and}\qquad u_{\epsilon}
\rightharpoonup u_{0}\quad\mathrm{in}%
\;\;H_{0}^{1}(\Omega)\quad\mathrm{for}
\;\epsilon\rightarrow0.
\end{equation*}
In order to obtain a strongly convergent sequence, the homogenization
  theory suggests to use  approximations with a correction,  namely,
\begin{equation*}
w_{\epsilon}^{1}\,({x})\,:=\,u_{0}({x})\,-\,\epsilon
\psi^{\epsilon}({x})\,N_{k}\,\Bigl(\frac{{x-x_{i}}%
}{\epsilon}\Bigr)\,\frac{\partial\,u_{0}({x})}{\partial x_{k}}%
\quad\forall{x}\in\Pi_{\mathbf{i}}^{\epsilon},\quad\forall\,{i},
\end{equation*}
where  $\,\psi^{\epsilon}:=\min\{1,\,\frac{1}%
{\epsilon}\;\mathrm{dist}({x},\,\partial\Omega)\}\,$
is a cutoff function. Then, optimal a priori convergence rates for the error
$u_{\epsilon}-w_{\epsilon}^{1}$ can be
proved  (e.g., see \cite{Bensoussan} Rem. 5.13, \cite{Cioranescu}, \cite{Jikov} p.28) if
$
u_{0}\,\in\,W^{2,\,\infty}(\overline{\Omega})$
and
$
\frac{\partial\,N_{k}}{\partial y_{j}}\in L^{\infty}(\widehat{\Pi}).
$
The resultant error estimate reads 
\begin{equation}
\label{1.8}
\Vert u_{\epsilon}-w_{\epsilon}^{1}\Vert_{H^{1}(\Omega)}\;\leq
\;\widetilde{c}\,\sqrt{\epsilon}. %
\end{equation}
Reconstruction of the flux
$A_{\epsilon}\,\nabla\,u_{\epsilon}$ with the same
convergence rate $\sqrt{\epsilon}$, requires solving
another  periodic problem for the operator
${\rm curl}\,\,A_{0}^{-1}{\rm curl}$.

In general,  above discussed correction techniques may be rather
costly and  the respective convergence  estimates usually require additional
assumptions concerning regularity of homogenized solutions. 
It uses solutions of boundary value problems on the cell
$\wh \Pi$ (e.g., (\ref{1.6})), which often can be found only approximately
and require analysis of effects generated
by approximation errors and their influence on the accuracy of
$\abig_0$, $u_0$, $w^1_\epsilon$, etc.
It should be also noted that the homogenization method provides accurate approximations
only for sufficiently small $\epsilon$ (this fact follows from the a priori estimate (\ref{1.8})). The question on
how to efficiently compute accurate approximations if $\epsilon$
is small but not "very small" remains open. One possible answer is
suggested below.

The approach considered in the paper is applicable to a much wider class of problems
than problems with periodically oscillating coefficients.
 It is valid for quasi-periodic structures 
(e.g. of the type presented in Fig. \ref{fig:3DPeriodStruct2}, 
 or periodic systems with defects, see \cite{VeKhorEwaldTuck:14}),
where homogenization theory cannot be used (see examples in \S5).
Structures of this type arise in various models in natural sciences and engineering applications
(see e.g. \cite{VeBoKh:Ewald:14,VeKhorCorePeriod:14} for applications in electronic structure calculations),
so that getting efficient approximations with guaranteed error bounds
is an important problem.


We  present a new computational strategy,
which is intended to compute efficiently approximate solutions
of boundary value problems with periodic and quasi-periodic structures in domains
composed of few tensor-product subdomains. This approach is based on 
the representation (approximation) of all vectors and matrices involved in the computational scheme 
in the so-called quantized TT (QTT) format \cite{KhQuant:09} such that all matrix-vector operations are 
implemented approximately via adaptive control the QTT rank parameters.
This method allows us to achieve the desired tolerance level 
regardless of the cell size $\epsilon$, i.e., it does not have limitations of the type (\ref{1.8}). 
Under certain assumptions on the 
tensor structure in coefficients and right-hand side the numerical cost can 
be estimated by $O(|\log^q \epsilon |)$, where the constant $q>0$ does not depend on 
$\epsilon$. 

In short, the main ideas behind this approach are as follows.
We use a simplified model with much simpler matrix $A_0$ as the basis
of the iteration algorithm (\ref{2.6}). If $A_\epsilon$ is defined by a perfectly regular
and highly oscillating structure, then
setting $A_0$ by (\ref{1.7}) is one possible option. However,  there are other options and the choice 
of $A_0$  is restricted only by the convergence conditions stated in
Theorem \ref{th2.1}. In other words, we can use  any matrix $A_0$ with simplified or averaged
coefficients that coarsely approximate the coefficients of $A_\epsilon$ if it 
provides contraction of the operator $T$ defined by (\ref{2.1}).
In more complicated cases, we can combine averaging, smoothing, and homogenization
in different parts of the domain in dependence of the structure and frequency of oscillations. 
The possible choice of $A_0$ is depicted in Fig. \ref{fig:3DPeriodStruct2} by red lines 
(see also examples in Fig. \ref{fig:Coef1DNperiod} and the discussion in Section  4.4).

The structure of $A_0$ defines the value contraction
factor $q$, which is explicitly
estimated a priori. Setting a  collection
of  simplified problems,  one can a priori
find the problem with minimal $q$ (this
amounts solving a simple optimization problem).
Next, on each step of the iteration algorithm we have guaranteed
two--sided a posteriori error estimates that control the distance
to the exact solution (see Section 3). 

It is important to outline that all steps of the iteration procedure
 and error estimates do not require
inversion of the matrix generated by the original (complicated)
differential operator. This matrix is involved
 only in  multiplication procedures, which can be performed very
fast and under in the small storage costs due to QTT tensor operations. This technique
is introduced in Section 4, which considers the low-rank quantized tensor 
approximations \cite{KhQuant:09}
 arising in the framework of our  tensor-based computational scheme. 
 We describe the QTT based preconditioned iteration
 and present the explicit QTT rank estimates for some particular classes  of equation coefficients in 1D. 
 We end up with numerical illustrations demonstrating
 the fast convergence of preconditioned iteration (geometric convergence rate), as well as 
 the logarithmic scaling of CPU computational time in the grid-size $N=2^L$.
Finally we note that solving differential equations with complicated and rapidly changing coefficients often
require special methods and approximations on distorted
meshes (see, e.g., \cite{Kuz1,Kuz2,KR2003} and references therein). We believe that  modifications of our
approach can be also helpful for these cases.

\section{The iteration method}\label{sec:PrecIter_Cont}
\subsection{Convergence}

We apply the general iteration method (see, e.g., \cite{GlLiTr})
in order to solve (\ref{1.1}) with the help of a simpler problem
generated by $A_0$. 
Let $v\in V_0$. The functional $\ell_v:V_0\rightarrow {\mathbb R}$
defined by the relation
$
\ell_v(w)=\abig_\epsilon(v,w)-(f,w)_\Omega
$
is a linear continuous functional. For any $\rho>0$ the problem: 
find $u\in V_0$ such that
\ben
\abig_0(u,w)
=\abig_0(v,w) - \rho \ell_v(w)\qquad\forall w\in V_0
\label {2.1}
\een
is well posed. 
Evidently, (\ref{2.1}) defines a linear bounded operator
$T:V_0\rightarrow V_0$ ($Tv=u$), which is
contractive
provided that $\rho$ is properly selected.
Indeed, select two arbitrary functions $v_1,v_2\in V_0$
and set
$u_1: = T(v_1)$, $u_2 := T(v_2)$. Let $e:=v_1-v_2$,
 and $\eta:=u_1-u_2$. Then 
 \begin{equation}
 \label{2.2}
 \abig_0(\eta,w)
=\abig_0(e,w) - \rho \ell_e(w)=\abig_0(e,w)-\rho \abig_\epsilon(e,w).
\end{equation}
 In view of (\ref{2.2}), the difference of images is subject
 to the relation
\begin{multline}
\label {2.3}
\|\eta\|^2_0:=\abig_0(\eta,\eta)
=\abig_0(e,\eta)-\rho
\abig_\epsilon(e,\eta)
=\IntO\left(A_0\nabla e-\rho\,
A_\epsilon \nabla e\right)\cdot \nabla \eta\,dx\\
\leq
\,\|\eta\|_{0}\left(\IntO (A_0\nabla e-\rho\,
A_\epsilon \nabla e)\cdot(\nabla e-\rho\,
A^{-1}_0A_\epsilon \nabla e)dx\right)^{1/2}.
\end{multline}
Hence,
\ben
\label{2.4}
\|\eta\|^2_0\leq\,\IntO (A_0-2\rho A_\epsilon+\rho^2
A_\epsilon A^{-1}_0A_\epsilon) \nabla e\cdot\nabla e\, dx
\leq\;q^2(\rho)\,\|e\|^2_0,
\een
where
 $q^2(\rho):=1-2\rho\lambda_1+\rho^2 c_\oplus$ and $c_\oplus$ is the constant in the estimate
\ben
\label{2.5}
A_\epsilon A^{-1}_0A_\epsilon\zeta\cdot\zeta\quad\leq
\quad c_\oplus A_0\zeta\cdot\zeta
\qquad\forall\zeta\in \Rd,\quad
x\in \Omega.
\een
If
$0<\rho<2\rho_*$, $\rho_*=\frac{\lambda_1}{c_\oplus}$ then the quantity
$q(\rho)$ is smaller than $1$
and (\ref{2.4}) yields the contraction estimate
$
\|\eta_0\|_{0}\,\leq\,q(\rho)\,\|e\|_{0}.
$
Well known results in the theory of fixed points (e.g., see \cite{Zeidler})
yield the following result.

\begin{theorem}
\label{th2.1}
For any $u_0\in V_0$
and $\rho\in (0,2\rho_*)$  the sequence $\{u_k\}$  of functions satisfying the relation
\ben
\label{2.6}
\abig_0(u_{k+1},w)=\abig_0(u_k,w)-\rho\,\ell_{u_k}(w)\qquad\forall w\in V_0,
\een
converges
to $u_\epsilon$ in $V$ and
$\|u_k-u_\epsilon\|_{0}
\leq q^k(\rho)\|u_0-u_\epsilon\|_{0}$
as $k\rightarrow+\infty$.
\end{theorem}
\subsection{Estimates of the contraction parameter}
For $\rho=\rho_*$
we find the contraction factor
\ben
\label {2.7}
q(\rho_*)=\sqrt{1-\frac{\lambda^2_1}{c_\oplus}}.
\een
Since
$
A^{-1}_0A_\epsilon\zeta\cdot A_\epsilon\zeta\;\leq
\; \frac{1}{\lambda^0_\ominus} A_\epsilon\zeta\cdot A_\epsilon\zeta\; \leq \;\frac{\lambda^\epsilon_\oplus}{\lambda^0_\ominus} 
A_\epsilon\zeta\cdot \zeta\;\leq\;
\frac{\lambda^\epsilon_\oplus}{\lambda^0_\ominus}\lambda_2
A_0\zeta\cdot\zeta,
$


we get a coarse estimate of the constant
$c_\oplus$ 
 and conclude that
\ben
\label {2.8}
q(\rho_*)\,\leq\,\sqrt{1-\frac{\lambda^2_1\lambda^0_\ominus}{\lambda_2\lambda^\epsilon_\oplus}}
=\sqrt{1-\left(\frac{\lambda^\epsilon_\ominus\lambda^0_\ominus}{\lambda^\epsilon_\oplus\lambda^0_\oplus}\right)^2}<1.
\een
However, in general the contraction factor can be better than in (\ref{2.8}).
Indeed, the value of $\rho$ in (\ref{2.4}) depends on the quantity
$\max\limits_{x\in \Omega} 
|{\mathbb B}^2(\rho,x,A_0,A_\epsilon)|$,
where
$${\mathbb B}(\rho,x,A_0,A_\epsilon):={\mathbb I}-\rho A^{-1}_0(x) A_\epsilon(x).
$$
If $A_0$ is given, then  the
best $\rho_*>0$ satisfies  the condition
\begin{equation}
\label{2.8a}
q(\rho_*,A_0,A_\epsilon)\,=\min\limits_{\rho\in {\mathbb R}_+}\max\limits_{x\in \Omega} |{\mathbb B}^2(\rho,x,A_0,A_\epsilon)|.
\end{equation}
In view of this fact, we obtain a criterium for selecting the best $A_0(x)$ among
a certain collection ${\mathcal A}$ of  "simple" (e.g., polynomial or piece-vise
constant) structures: find $\rho_*>0$ and $A^*_0\in\mathcal A$ such that
\begin{equation}
\label{2.9}
q(\rho_*,A^*_0,A_\epsilon):=\max\limits_{x\in \Omega} |{\mathbb B}^2(\rho_*,x,A^*_0,A_\epsilon)|\,=
\min\limits_{\,A_0\in {\mathcal A}}
\min\limits_{\rho\in {\mathbb R}_+}\, \max\limits_{x\in \Omega} |{\mathbb B}^2(\rho,x,A_0,A_\epsilon)|.
\end{equation}
In the right hand side of (\ref{2.9}) we have a matrix optimization problem. Solving it
leads to the best simplified matrix $A_0$ and optimal $\rho_*$. It is worth outlining that
this relatively simple problem does not require solutions of some auxiliary boundary
value problems and can be done before
iterative  computations
based on (\ref{2.6}).

{\em Particular case.} Let $A_0=a_0{\mathbb I}$ and $A_\epsilon=a_\epsilon{\mathbb I}$. In this case,
 (\ref{2.8}) means that $\rho_*$ should be selected such that
\begin{equation}
\label{2.10}
\max\limits_{x\in \Omega}
|1-\rho_* {\mathds h}(x)|=
\min_\rho \max\limits_{x\in \Omega}|1-\rho {\mathds h}(x)|,\qquad {\mathds h}(x)=\frac{a_\epsilon}{a_0}>0.
\end{equation}
If $d=1$, then it is not difficult to show that
$$
\max\limits_x|1-\rho {\mathds h}(x)|=\max\{|1-\rho \underline {\mathds h}|,|1-\rho \overline {\mathds h}|\},
$$
where $\overline {\mathds h}:=\max\limits_{x\in \Omega}{\mathds h}(x)$ and
$\underline {\mathds h}:=\min\limits_{x\in \Omega}{\mathds h}(x)$.
We find that $\rho_*=\frac{2}{\underline {\mathds h}+\overline {\mathds h}}$ and
 $$
 q(\rho_*,a_0,a_\epsilon)=\frac{\overline {\mathds h}-\underline {\mathds h}}{\underline {\mathds h}+\overline {\mathds h}}<1.
 $$
For $a_0=const$, we obtain
$q(\rho_*,a_0,a_\epsilon)=\displaystyle{
\frac{
\overline{a}_\epsilon-
\underline{a}_\epsilon}
{\overline{a}_\epsilon+
\underline{a}_\epsilon}}$.
 
If ${\mathcal A}$ is a set of possible 
simplified coefficients, then the optimization problem (\ref{2.9}) is reduced to finding $a^*_0\in {\mathcal A}$
that solves the problem
\begin{equation}
\label{2.11}
\min\limits_{a_0\in \mathcal A}\left(\frac{\max\limits_{x\in \Omega}
\frac{a_\epsilon}{a_0}-\min\limits_{x\in \Omega}\frac{a_\epsilon}{a_0}}
{\max\limits_{x\in \Omega}\frac{a_\epsilon}{a_0}+
\min\limits_{x\in \Omega}\frac{a_\epsilon}{a_0}}\right).
\end{equation}

Let $a_\epsilon$ be a function
oscillating around a certain "mean" function $a_0$. Assume that
maximal relative deviation of $a_0$
from $a_\epsilon$ does not exceed $\delta$, i.e.,
$
\displaystyle\frac{a_\epsilon}{a_0}\in [1-\mu_\ominus,1+\mu_\oplus]$,
$\mu_\ominus\in (0,1)$, $\mu_\oplus\geq 0$. The paprameter
 $\delta:=\mu_\ominus+\mu_\oplus$ charackterises the scale of deviations.
Since
$
\max\limits_{x\in \Omega}
\frac{a_\epsilon}{a_0}\leq 1+\mu_\oplus,
$
and
$\min\limits_{x\in \Omega}
\frac{a_\epsilon}{a_0}=1-\mu_\ominus$,
we find that 
$$
 q(\rho_*,a_0,a_\epsilon)=\frac{\delta}{2+\mu_\oplus-\mu_\ominus}=
 1-\frac{1-\mu_\ominus}
 {\delta/2+1-\mu_\ominus}.
$$
This formula shows that the method can be very efficient if $\delta$
is small, i.e., if $a_\epsilon$
oscillates around $a_0$ with a relatively small amplitude.

\subsection{Discrete setting}

Theorem (\ref{th2.1}) is applicable to the case where the problem
is solved on a finite dimensional subspace $V_{0h}\in V_0$ associated with a mesh $\cT_h$.
Let $u_{0,h}$ solve the problem
\be
\abig_0(u_{0,h},w_h)=(f,w_h)_\Omega\qquad\forall w_h\in V_{0h}
\ee 
and the functions $\{u_{k+1,h}\}\in V_{0h}$, $k=0,1,2,...,$ satisfy
\ben
\label{2.13}
\abig_0(u_{k+1,h},w_h)=\abig_0(u_{k,h},w_h)-\rho\,\ell_{u_{k,h}}(w_h)\qquad\forall w_h\in V_{0h}.
\een
 By repeating the same arguments as before, we conclude that $\{u_{k,h}\}$ tends to the fixed point $u_{\epsilon,h}$ 
 of (\ref{2.13}) provided that $\rho<2\rho_*$. Obviously,
$u_{\epsilon,h}$ satisfies the relation
 \ben
 \label{2.14}
0= \ell_{u_{k,h}}(w_h):=
\abig_\epsilon
(u_{\epsilon,h},\nabla w_h)-(f,w_h)_\Omega\qquad\forall w_h\in V_{0h}.
 \een
 We see that $u_{\epsilon,h}$ is the Galerkin
approximation of $u_\epsilon$ on $V_{0h}$
and obtain the a priori error estimate
\ben
\label{2.15}
\|u_{k,h}-u_{\epsilon,h}\|_{0}
\leq q^k(\rho)\|u_{0,h}-u_{\epsilon,h}\|,
\een
which shows that approximations converge to $u_{\epsilon, h}$ with the geometric rate.
>From the practical point of view it is more important
to have an estimate of $\|u_{k,h}-u_\epsilon\|_0$.
In the next section we will obtain such estimates.

Now we discuss matrix equations that follow from
(\ref{2.13}) and compare them with the equations
generated by  a "straightforward" approach applied
to (\ref{1.3}).
Let $\{\phi_i\}$, $i=1,2,...N$ be a system of linearly independent trial functions and
$V_h={\rm Span}\{\phi_i\}$. Define the vector
$\f:=\{f_i\}$, $f_i= (f, \phi_i)_\Omega$ 
and two matrixes
${\Bbb A}_\epsilon:=
 \{a_\epsilon (\phi_i,\phi_j)\}$ and
${\Bbb A}_0:=\{a_0 (\phi_i,\phi_j)\}.
$
In particular, we can use piecewise affine
basis functions in $V_{0h}\subset H^1_0(\Omega)$ associated with 
the uniform tensor-product Cartesian grid.
Denoting the fine grid size by $h=1/(N+1)$, where $N$ is the number
of grid points in each spatial direction, the total problem size is estimated by $N^d$.
For ease of exposition we assume that each scaled unit cell of univariate size 
$O(\epsilon)$ includes equal number $n_0$ of grid points.

Direct computation of $u_{\epsilon,h}$ requires solving the algebraic problem
\begin{equation}\label{2.16}
 \mathbb{A}_\epsilon {\bf v}_\epsilon = {\bf f}, \qquad \mathbb{A}_\epsilon \in \mathbb{R}^{N^d \times N^d},
 \quad  {\bf f} \in \mathbb{R}^{N^d},
\end{equation}
with a sparse stiffness matrix $\mathbb{A}_\epsilon$. 
Here ${\bf v}_\epsilon  \in \mathbb{R}^{N^d}$
is the vector of nodal values 
that define $u_{\epsilon,h}$.
The main bottleneck of the above computational scheme is due to the 
matrix size in the Galerkin system (\ref{2.16}).
Indeed, the univariate mesh parameter $N$ 
is of the order of $N=O(\frac{n_0}{\epsilon})$, where $n_0$ is the mesh parameter that
ensures the sufficient resolution of all data in the  cell of length $\epsilon$. In general, this parameter
also depends on another structural parameter $\kappa$.
Hence accurate approximations require huge values of $N$ so that the numerical complexity of the direct solver for the system (\ref{2.16})
scales polynomially in the frequency parameter $(1/\epsilon)^d$.
Homogenization methods introduce a model
simplification providing indirect $O(\sqrt{\epsilon})$-approximation to the 
solution of (\ref{2.16}). This method avoids inversion
of $\mathbb{A}_\epsilon$ and leads to the problem (\ref{1.7}), which can be solved on subspaces of much lower dimension.

The iteration scheme (\ref{2.13})
suggests another way to avoid inversion of $A_\epsilon$.
The basic iteration algorithm  on the full 
finite element space (exact arithmetics)
starts with $\v_0={\Bbb A}^{-1}_0\f$ and computes
$\v_{k+1}$, $k=0,1,2, ...$ by solving the problem
\begin{equation}
 \label{2.17}
 {\Bbb A}_0 \v_{k+1}={\Bbb A}_0 \v_k-\rho({\Bbb A}_\epsilon \v_k - \f).
\end{equation}
We can rewrite (\ref{2.17}) in the form
$\v_{k+1}- \v_k=
\rho({\Bbb A}^{-1}_0 \f-{\Bbb A}^{-1}_0
{\Bbb A}_\epsilon \v_k)$, 
which shows that (\ref{2.17}) is equivalent to the iteration method applied to the preconditionered system
\ben
\label{2.18}
{\mathbb A}^{-1}_0{\Bbb A}_\epsilon {\bf u}_\epsilon=\v_0.
\een

It is worth awaiting that in many
cases (at least for periodic or almost periodic coefficients with small $\epsilon$) the above selection 
of  $\v_0$ with ($\mathbb{A}_0$ generated by the homogenized
problem (\ref{1.7}) or another suitable simplified matrix) will  provide a good starting approximation to the procedure (\ref{2.17}).
 This fact was  indeed confirmed in various numerical tests which show that such a constructed $\v_0$
 is a good initial guess  for the iteration method.

 The iteration (\ref{2.17}) involves only one
operation with $\mathbb{A}_\epsilon$: multiplication
by the vector $\v_k$. If entries of the matrix $\mathbb{A}_\epsilon$ are 
generated by oscillating functions having low rank approximation in the so-called quantized 
tensor representations (QTT) \cite{KhQuant:09}, then (\ref{2.17})
can be solved fairly easily by QTT-based tensor type methods applied to 
the properly transformed linear system (see Section 4). 
The approximate tensor arithmetics makes performing the iterations inexpensive. 
Thus, we obtain a new computational approach for a rather wide class
of problems with periodic and quasi--periodic coefficients
that allows to solve equation (\ref{2.16}) by iteration (\ref{2.17}) 
with the required precision at the cost that scales only logarithmically in $\epsilon$.


\section{Error control}
\label{ssec:Error}

For the control of approximation errors we use a posteriori estimates
of the functional type (see \cite{Re2000,ReGruyter} and the references
therein). They provide guaranteed and fully computable
bounds of errors for any conforming
approximation within the framework of a unified procedure, which does not
require special features of approximations (e.g.,
exact satisfaction of the Galerkin orthogonality condition) or
special features of the exact solution (e.g., extra regularity).
Such estimates are
robust and convenient  for problems with complicated
coefficients (see, e.g., \cite{ReSaSa1}), where the above mentioned conditions are difficult to guaranty.
We recall that for any approximation $v\in V$ of the problem
(\ref{1.1}) we have the following  estimate
of the error $e:=u_\epsilon-v$
\begin{equation}
 \label{3.1}
\| e\|_{\epsilon}\leq\,
\|A_\epsilon \nabla v-y\|_{\epsilon}+
C_\Omega\|{\rm div}y+f\|:=M_\oplus(v,y),
\end{equation}
where $\|e\|^2_{\epsilon}=\abig_\epsilon(e,e)$,
$y$ is an arbitrary vector  function in $H(\Omega,{\rm div})$
and $C_\Omega$ is the Friedrichs constant. If
$
\Omega\subset\left\{x\in {\mathbb R}^d\,\mid\,a_s<x<b_s,\;
b_s-a_s=l_s,\;s=1,2,..,d\right\},
$
then
$C_\Omega=\frac{1}{\kappa\pi}$, $
\kappa^2=\sum\limits^d_{s=1}\frac{1}{l^2_s}.
$
In particular,
$C_\Omega=\frac{1}{d^{1/2}\pi}$  for $\Omega=(0,1)^d$.

The functional $M_\oplus(v,y):H^1_0(\Omega)\times H(\Omega,{\rm div})\rightarrow {\mathbb R}_+$  is an error majorant. 
Properties of such error majorants are well studied
(see, e.g., \cite{ReGruyter} ). We know that
$M_\oplus(v,y)$ vanishes
if and only if $v=u_\epsilon$ and $y=p_\epsilon:=A_\epsilon\nabla u_\epsilon$.
Moreover, for any $v\in H^1_0(\Omega)$ the functional
$M_\oplus(v,p_\epsilon)$ coincides with the error
and the integrand of $M_\oplus(v,p_\epsilon)$ shows the distribution of local errors.
 Numerous
tests performed for different boundary value problems have
confirmed practical efficiency of this and other error majorants
derived for various problems.  It was shown
 that  $M_\oplus$ is a guaranteed and efficient majorant of the global 
error and good indicator of local 
errors if the exact flux is replaced by a certain numerical reconstruction $p_h$ (in our case 
instead of $p_\epsilon$ we use
$p_{h,\epsilon}$). 
There are many different ways to obtain suitable reconstructions (e.g., see \cite{MaNeRe} for a systematic discussion of
computational aspects of this error estimation
method). 

Error majorants can be efficiently used for the evaluation of modeling errors (see \cite{ReSaSm,ReSaSa2}).
In particular, if we set $y=A_0\nabla u_0$, then (\ref{3.1})
implies a simple estimate of the modeling error caused
by simplification of coefficients:
$
\| e\|_{\epsilon}\leq\,
\|(A_\epsilon-A_0) \nabla v\|_{\epsilon}.
$

However, for problems with highly oscillating coefficients the general
majorant (\ref{3.1}) has a substantial drawback: it contains a norm 
generated by $A^{-1}_\epsilon$. In our analysis we try to avoid all the operations
related to this most complicated matrix except multiplication by a vector (which can be
performed by tensor type methods). Hence, the goal is to
modify general  a posteriori estimates in accordance with this principle.
Consider one step of the iteration method, where the function $v\in V_0$ generates $u=Tv\in V_0$.
For a contractive mapping $T$  we have two sided error
estimates (see \cite{Ostrowski,ReGruyter,Zeidler}), which imply

\ben
\label{3.2}
\|v-u_\epsilon\|_0\,\in
\left\{\frac{1}{1+q(\rho)},\,\frac{1}{1-q(\rho)}\right\}\|u-v\|_0
\een
Here $v$ is known, but $u=Tv$ is generally unknown and we need to use some approximation
$\wt u$ instead. 
It is easy to see that
\ben
\label{3.3}
\|u-v\|_0\,\leq\, \|\eta\|_0+\|u-\wt u\|_0,\\
\label{3.4}
\|u-v\|_0\,\geq\, \|\eta\|_0-\|u-\wt u\|_0,
\een
where the function $\eta:=\wt u-v$ is known. The estimates (\ref{3.3})
and (\ref{3.4}) would be fully computable if we find a computable
majorant of the norm $\|\wt u-u\|_0$.

We note that
\ben
\label{3.5}
\abig_0(u,w)=\abig_0(v,w)-\rho\,
\IntO( A_\epsilon \nabla v\cdot \nabla w-fw)\qquad\forall w\in V_0.
\een
Then for any $y\in H(\Omega,{\rm div})$ and $w\in V_0$, we have
\ben
\label{3.6}
&\abig_0(u-\wt u,w)&=\abig_0(v-\wt u,w)-\rho\,
\IntO( A_\epsilon \nabla v\cdot \nabla w-fw)\\
&&\IntO\left(A_0\nabla (v-\wt u)\cdot \nabla w-\rho(A_\epsilon \nabla v\cdot \nabla w-fw\right))dx
\nonumber\\
&&\IntO(A_0\nabla (v-\wt u)-\rho A_\epsilon \nabla v+y)\cdot\nabla w+
(\rho f+{\rm div} y)w)dx.
\nonumber
\een
Using the notation $\tau:=y-\rho A_\epsilon \nabla v$, we obtain
\begin{multline}
\label{3.7}
\|u-\wt u\|_0\leq\,
\left(\IntO (
A_0\nabla \eta\cdot\nabla \eta+
A^{-1}_0\tau\cdot  \tau-2\nabla\eta\cdot\tau )dx\right)^{1/2}\\
+\frac{C_F}{\lambda^0_\ominus}
\|{\rm div} y+\rho f\|=:\cM_\oplus(\wt u,v,y).
\end{multline}
Note that
\be
\inf\limits_{y\in H(\Omega,{\rm div}}
\cM_\oplus(\wt u,v,y)=\|u-\wt u\|_0.
\ee
Indeed, set $y=A_0\nabla(u-v)+\rho A_\epsilon\nabla v$.
Then, $\tau=A_0\nabla(u-v)$, ${\rm div} y+\rho f=0$,
and the first term of the majorant is equal to $\|u-\wt u\|_0$.
Hence, the estimate has no gap.

Now we apply these relations to the step $k$ of (\ref{2.6}).
Set $v=u_{k,h}$ (approximation computed at step $k$
using a mesh $\cT^k_h$). Then $u=Tu_{k,h}$
is the exact solution of (\ref{3.5}), which we do not know.
Instead we have a function $\wt u=u_{k+1,h}$ computed 
on the mesh $\cT^{k+1}_h$ (it may coincide with the previous mesh $\cT^{k}_h$ or be a new one constructed by, e.g., a refinement procedure). 
The function $\eta=\eta_{k+1}:=u_{k+1,h}-u_{k,h}$ is known.
By (\ref{3.3}) and (\ref{3.4}) we obtain
\ben
\label{3.8}
&&\|u_{k,h}-u_\epsilon\|_0
\leq\, \frac{1}{1-q}\left(\|\eta_{k+1}\|_0+
\cM(u_{k+1,h},u_{k,h},y)\right),\\
\label{3.9}
&&\|u_{k,h}-u_\epsilon\|_0\geq\,
\frac{1}{1+q}\left(\|\eta_{k+1}\|_0-\cM(u_{k+1,h},u_{k,h},y)\right).
\een
Here $y$ is any vector function in $H(\Omega,{\rm div})$.
Certainly, getting minimal values of the majorant
require  a suitable numerical reconstruction 
of the exact flux of the problem (\ref{2.6}), which is
$
q_k=A_0\nabla u_{k+1}-\sigma_k$,
where $\sigma_k=(A_0-\rho A_\epsilon)\nabla u_{k,h}$
is known. Since the coefficients of $A_0$  are regular
and do not oscillate, reconstructions of such a flux
can be done by well known methods (see, e.g., \cite{MaNeRe,ReGruyter} and the literature cited in these
books).
We do not discuss this question in detail because it
will be the matter of a special publication focused
on multidimensional problems. 

If $d=1$ 
(related
to the numerical tests below), then the
flux is easy to reconstruct.
In this case the problem (\ref{1.1}) is $(a_\epsilon u^\prime_\epsilon)^\prime+f=0$ and the simplified problem
is $(a_0 u^\prime_0)^\prime+f=0$.
We set $y=\rho(-g(x)+c)$ (where $g(x)=\int^x_0fdx$) and define the constant $c$ by minimizing the first term of ${\mathcal M}(u_{k+1,h},u_{k,h},y)$.
We have $\tau=\rho(c-g(x)- a_\epsilon u^\prime_k)$
and need to find $c$ minimizing the quantity
\be
\int\limits^1_0 (a_0|\eta^\prime_{k+1}|^2 +a^{-1}_0\rho^2(c-g(x)- a_\epsilon u^\prime_{k,h})^2-2\rho\eta^\prime_{k+1} (c-g(x)- a_\epsilon u^\prime_{k,h})dx.
\ee
 Since $\int^1_0 \eta^\prime_{k+1}dx=0$, the problem is reduced to
 minimization of the second term
and the best
 $c$ satisfies the equation
$
\int\limits^1_0 a^{-1}_0
(c-g(x)-a_\epsilon u^\prime_{k,h})dx=0.
$
Hence
$
c=c_k:=\left(\int\limits^1_0 a^{-1}_0(g(x)+a_\epsilon u^\prime_{k,h})dx\right)\left(\int\limits^1_0 a^{-1}_0\,dx\right)^{-1}
$
and
\begin{multline*}
\cM^2(u_{k+1,h},u_{k,h},y)=\\
=\int\limits^1_0 (a_0|\eta^\prime_{k+1}|^2 +a^{-1}_0\rho^2(c_k-g(x)- a_\epsilon u^\prime_{k,h})^2+2\rho\eta^\prime_{k+1} (g(x)+a_\epsilon u^\prime_{k,h}))dx.
\end{multline*}
We see that the majorant is fully computable.
Moreover, if the sequence $\{u_{k,h}\}$ converges to $u_\epsilon$ in $V$, then $\|\eta_{k+1}\|_0\rightarrow 0$. Hence the first and the last terms
of the above integral tend to zero. Also,
$g(x)+a_\epsilon u^\prime_{k,h}\rightarrow g(x)+a_\epsilon u^\prime_\epsilon=:c_*$ in $L^2$. 
Note that
\be
c_k=\frac{\int\limits^1_0 a^{-1}_0(g(x)+a_\epsilon u^\prime_{k,h})dx}{\int\limits^1_0 a^{-1}_0\,dx}=
c_*+\frac{\int\limits^1_0 a^{-1}_0a_\epsilon (u^\prime_{k,h}-u^\prime_\epsilon)dx}{\int\limits^1_0 a^{-1}_0\,dx}
\ee
Therefore,
$
|c_k-c_*|
\leq \mu_{0,\epsilon}
\|u^\prime_{k,h}-u^\prime_\epsilon\|_0
$, where
$\mu_{0,\epsilon}=
\left(\int\limits^1_0 a^{-1}_0a^2_\epsilon dx\right)^{1/2} \left(\int\limits^1_0 a^{-1}_0\,dx\right)^{-1}.
$
does not depend on $k$. Therefore, the second term also tends to zero and we conclude
that  the majorant (\ref{3.8}) tends to zero. Also, it is possible to show that the majorant is equivalent
to the error.

\section{Tensor-based preconditioned iterative scheme }\label{sec:Tens_Precond_Iter}

The main concept of our tensor-based approach is the direct iterative 
solution of the initial large algebraic system (\ref{2.16})
in the form of preconditioned iteration (\ref{2.17})
 by using low-parametric data formats,  exploiting 
certain redundancy in the grid-based representation of matrices and vectors involved.
This is realized, first, by transformation of the ``low-dimensional'' 
FEM-Galerkin equation to the equivalent system posed in the high dimensional
quantized  tensor space, and then by solving this system iteratively using the low-rank
QTT tensor approximation \cite{KhQuant:09} to the Galerkin stiffness matrix, 
the preconditioner and all vectors involved.
This allows to compute the numerical approximation to the exact solution 
discretized on a fine grid
up to the chosen precision $\delta>0$, adapted to the mesh-resolution but independent
on the frequency parameter $1/\epsilon$. 

The approach is well adapted 
to fast  QTT-based tensor approximation method, what is natural
to await because the QTT tensors fit well the intrinsic features 
of FEM discretizations to functions and operators generated via periodic 
and quasi--periodic geometric structures 
\cite{KhQuant:09,VeBoKh:Ewald:14,VeKhorEwaldTuck:14,KhorVeit:14}.
The numerical cost of the rank-structured 
iteration can be bounded by $O(|\log \epsilon|^q)$ provided that rank parameters
remain small.

\subsection{QTT tensor representation of function related vectors  and matrices} 
\label{ssec:QQTform}

In this section we present a brief overview of QTT tensor approximation 
method \cite{KhQuant:09} which is the base for the construction of the presented 
tensor-based computational scheme.
We refer to surveys on commonly used  low-rank representations
of discrete functions and operators \cite{Kolda,KhorSurv:10,Scholl:11}.

The QTT-type approximation of an $N$-vector with $N=q^L$, $L\in \mathbb{N}$ (usually $q=2$) 
allows to reduce the asymptotic storage cost to $O(\log N)$ \cite{KhQuant:09}.
The QTT rank decomposition applies to a tensor
obtained by the $q$-adic folding (reshaping) of the target long vector to an $L$-dimensional
$q\times ... \times q$ data array considered in the $L$-dimensional quantized tensor space.
As the basic result, in \cite{KhQuant:09} it was shown that for a large class 
of function related vectors (tensors)
such a procedure allows the low-rank representation of their quantized $L$-dimensional image, 
thus reducing the representation complexity to the logarithmic scale $O(\log N)$.

In particular, a vector
$
{\bf x}=[x(i)]_{i=1}^N\in \mathbb{R}^{N},
$
is reshaped to its quantics image in
$
\mathbb{Q}_{q,L}= \bigotimes_{j=1}^{L}\mathbb{K}^{q},
 \; \mathbb{K}\in \{\mathbb{R},\mathbb{C}\},
$
by $q$-adic folding,
\[
\mathcal{F}_{q,L}: {\bf x} \to \textsf{\textbf{Y}}
=[y({\bf j})]\in \mathbb{Q}_{q,L}, \quad {\bf j}=\{j_{1},...,j_{L}\},
\quad \mbox{with} \quad j_{\nu}\in \{1,2,...,q\}, 
\]
where for fixed $i$, we have $y({\bf j}):= x(i)$, and $j_\nu=j_\nu(i)$, $\nu=1,...,L$,
is defined via $q$-coding,
$
j_\nu - 1= C_{-1+\nu},
$
such that the coefficients $C_{-1+\nu} $ are found from the
$q$-adic representation of $i-1$,
\[
i-1 =  C_{0} +C_{1} q^{1} + \cdots + C_{L-1} q^{L-1}\equiv
\sum\limits_{\nu=1}^L (j_{\nu}-1) q^{\nu-1}.
\]

Suppose that the quantized image for certain  $N$-vector
(i.e. an element of $L$-dimensional quantized tensor space $\mathbb{Q}_{q,L}$ with $L=\log_q N$)
can be effectively represented (approximated) in the low-rank canonical or TT format.
For given QTT-rank parameters $\{r_k\}$ ($k=1,...,L-1$) 
the number of  representation parameters in the QTT approximation  can be estimated by
$$
q r^2 \log_q N \ll N, \quad \mbox{where}\quad r_k \leq r, \quad k=1,...,L-1,
$$
providing $\log$-volume scaling in the size of initial vector, $N$.
The optimal choice of the base $q$ is shown to be $q=2$ or $q=3$ \cite{KhQuant:09}, however
the numerical realizations are usually implemented  by using binary coding, i.e. for $q=2$.
For $d\geq  2$ the construction is similar \cite{KhQuant:09}.

The favorable features of the QTT approximation method are due to the
perfect low rank decompositions discovered for 
the wide class of function-related tensors \cite{KhQuant:09}: 
\begin{proposition} \label{prop:KhQuant_Appr:09}  (\cite{KhQuant:09})
 Let vector ${\bf x}\in \mathbb{C}^N$, $N=2^L$, be obtained by sampling a continuous 
 function $f\in C[0,1]$ over the uniform grid
 of size $N$. Then the following QTT-rank estimates are valid independent
 on the vector size $N$:\\
 (A) $r=1$ for complex exponentials, $f(x)=e^{i \omega x}$, $\omega \in \mathbb{R}$.\\
 (B) $r=2$ for trigonometric functions, $f(x)=\sin {\omega x}, f(x)=\cos {\omega x}$, 
 $\omega \in \mathbb{R}$. \\
 (C) $r\leq m+1$ for polynomials of degree $m$. \\
 (D)
For a function $f$ with the QTT-rank $r_0$ modulated by another function $g$ with 
the QTT-rank $r$ (say, step-type function, plain wave, polynomial) 
the QTT rank of a product $f\, g$ is bounded by a multiple of $r$ and $r_0$.\\
 (E) QTT rank for the periodic amplification of a 
reference function on a unit cell to a rectangular lattice is of the same order as that 
for the reference function \cite{VeBoKh:Ewald:14}.
\end{proposition}

Concerning the matrix case,
it was found in \cite{Osel-TT-LOG:09} by numerical tests
that in some cases the dyadic reshaping of an $N \times N$ matrix with $N=2^L$ may lead to a small
TT-rank of the resultant  matrix rearranged to the tensor form.
The explicit low-rank QTT representations for a class of discrete multidimensional
matrices mapping the space $\mathbb{Q}_{2,L}$ into itself  were proven in \cite{KazKhor_1:10}, 
see also \cite{KhorSurv:10}.

In our applications the concept for the construction of fast numerical methods
is based on the $\epsilon$-independent low-rank QTT approximation of all function-related 
vectors and operator-related matrices involved in the computational scheme. 
This allows to reduce the numerical 
costs to the logarithmic scale in the grid-size, i.e. of order $O(|\log \epsilon|^q)$.
The critical issue concerns with QTT approximation of the FEM-Galerkin matrix 
$\mathbb{A}_\epsilon$ generated by 
the highly oscillating or/and jumping diffusion coefficient, and of the 
respective ``homogenized'' preconditioner $\mathbb{A}_0$.

\subsection{The Galerkin FEM scheme}

For further discussion we choose  the Galerkin FEM with $N$ piecewise-linear hat 
functions $\left\{\phi_i\right\}$
in the physical domain $\Omega =[0,1]$, constructed on a fine uniform grid with 
step size $h=1/(N+1)$, 
which is a small fraction of $\epsilon$,
and nodes $x_{i}=hi$, $i = 1,\ldots,N$. For ease of notation we further 
denote $a_\epsilon(x)=A_\epsilon(x)$, then
the entries of the exact stiffness matrix $\mathbb{A}[a_\epsilon]$  read
\begin{equation}\label{eqn:GalekMatr}
\left(\mathbb{A}[a_\epsilon]\right)_{i,i'} = \left(a_\epsilon(x)\nabla\phi_i(x), 
\nabla\phi_{i'}(x)\right)_{L_2(D)}, \quad i,i'=1,\ldots,N.
\end{equation}

To simplify the approximation procedure, we may assume that the coefficient 
remains constant at each
spatial interval $[x_{i-1}, x_i]$, which corresponds to the evaluation of the scalar product
above via the midpoint quadrature rule.
It is known that this quadrature yields the approximation order $\Order{h^2}$, the same as
the piecewise-linear discretization of the solution.

We introduce the coefficient vector
${\bf a}=[a_{i}]\in \mathbb{R}^N$, $a_{i} = a_\epsilon(x_{i-1/2})$, $i=1,\ldots,N$, then
the resulting tridiagonal matrix takes the form,
\begin{equation}
\mathbb{A}[{\bf a}] = \dfrac{1}{h}\left[
\begin{matrix}
a_1+a_2 & -a_2 \\
-a_2 & a_2+a_3 & -a_3 \\
& \ddots & \ddots & \ddots \\
& & -a_{N-1} & a_{N-1}+a_{N} & -a_{N} \\
& & & -a_{N} & 2a_{N}
\end{matrix}\right].
\label{eqn:stiff_matrix_1point}
\end{equation}

\subsection{QTT tensor representation of the system matrix}\label{ssec:QQTformMatr}

In this section we discuss the low-rank QTT representations of the Galerkin matrices
approximating elliptic operators with variable coefficients in 1D.
In particular, we construct the QTT representation of arising three-diagonal
stiffness matrices by using the related results in \cite{DoKazKh_1DSPDE:12}.
The approach can be extended to $d$-dimensional equations defined  
on the lattice-type geometries.

We let $N=2^L$ and represent the indices in physical space in the binary coding
$$
i = \overline{i_1,\ldots,i_{L}}, \quad \mbox{where}\quad L=\log_2 N,
$$
to consider the QTT decomposition of vectors and matrices involved in the discrete problem.
Then the coefficient vector ${\bf a}=[a_i]$ in \eqref{eqn:stiff_matrix_1point} can be
represented in the rank-${\bf r}$ QTT form as an $L$-dimensional tensor, where 
${\bf r}=(r_1,...,r_{L-1})$:
\begin{equation}
\begin{split}
{a}_i  & = \sum\limits_{k_1,\ldots,k_{L-1}=1}^{r_1,\ldots,r_{L-1}}
a^{(1)}_{k_1}(i_{1}) \cdots a^{(L)}_{k_{L-1}}(i_{L}).
\end{split}
\label{eqn:a_qtt_sep}
\end{equation}

The QTT cores of the matrix $\mathbb{A}[{\bf a}]$ can be written similarly, 
recalling that the vector is turned to the diagonal matrix without changing the TT ranks.
Specifically,
the matrix $\mathbb{A}[{\bf a}]$ can be brought into the QTT format by using the 
shift matrices and their explicit rank-2 QTT representation.
Let us denote by $\mathbb{S}=[s_{i,i'}] \in \mathbb{R}^{N\times N}$ the upper shift matrix
given by
$$
s_{i,i'} = \left\{\begin{array}{ll} 1, & i'=i+1, \\ 0, & \mbox{else,} \end{array}\right. ,
$$
and notice that
this matrix has exact rank-$2$ QTT representation \cite{KazKhor_1:10}. Then it holds
\begin{equation}\label{eqn:Gamma_x}
\mathbb{A}\left[{\bf a} \right] = \mathbb{S} \diag({\bf a}) + \diag({\bf a}+ \mathbb{S} {\bf a})
 + \diag({\bf a}) \mathbb{S}^\top,
\end{equation}
with the maximal QTT ranks estimate
$
r \left(\mathbb{A}\left[ {\bf a} \right] \right) \le 7 r({\bf a}),
$
controlled  by the QTT rank of the coefficients vector ${\bf a}$.
The representativity of the stiffness matrix in low-rank formats can be summarized 
as follows \cite{DoKazKh_1DSPDE:12}.
\begin{theorem}\label{thm:QTTrMatr}
Let the problem
is discretized by the matrix (\ref{eqn:stiff_matrix_1point}) with the use of the
Galerkin-FEM method on the uniform grid. Suppose that the diffusion coefficient
vector ${\bf a}$ is given in a form of a QTT decomposition \eqref{eqn:a_qtt_sep}
with the following QTT-rank bounds
$r_p \le R$, for $p=1,\ldots,L-1$.
Then the  QTT ranks of the matrix $\mathbb{A}[{\bf a}]$ can be bounded by $7 R$.
\end{theorem}

Theorem \ref{thm:QTTrMatr} ensures that the arising Galerkin system of linear equations
(\ref{2.16})
can be efficiently solved iteratively as (\ref{2.17}) by using the low-rank QTT representation of 
each of entities ${\bf a}, {\bf u}_\epsilon,$ and ${\bf f}$,
provided that the spectrally close preconditioner, $\mathbb{A}_0$, to the stiffness 
matrix $\mathbb{A}_\epsilon\left[{\bf a}\right]$ is constructed. 

\subsection{Main assumptions and the construction of preconditioner}\label{ssec:Prec_Assumpt}

In the following numerical tests we chose the homogenized elliptic operator 
defined by the coefficient $a_0(x)$ as the preconditioner.

In the case of exotic coefficients, the construction of ``homogenized'' 
coefficient may depend on the shape of the initial equation coefficient $a_\epsilon(x)$.
We define the preconditioning operator $\mathbb{A}_0$ via the generalized 
averaging procedure
$$
\widetilde{a}_0(x)=\frac{1}{2}(a^+(x) + a^-(x)),
$$ 
where $a^+(x)$ and $a^-(x)$ are chosen as {\it majorants and minorants} 
of $a_\epsilon(x)$, respectively.
Examples of such a construction is given in Figure \ref{fig:3DPeriodStruct2}, left and middle,
and in Figure \ref{fig:Coef1DNperiod}, where the coefficient $\widetilde{a}_0(x)$
is colored in red.
The estimate on the condition number of the preconditioned matrix is 
given by the following simple lemma.
\begin{lemma}
Define $q(x):=(a^+(x) - \widetilde{a}_0(x))/\widetilde{a}_0(x)$, then
 the condition number of the preconditioned matrix
$\mathbb{A}_0^{-1} \mathbb{A}_\epsilon$ is bounded by
\[
 cond \{\mathbb{A}_0^{-1}\mathbb{A}_\epsilon \} \leq C \max\frac{1+q(x)}{1-q(x)}.
\]
\end{lemma}
The main assumptions for applicability of the presented tensor method are the following:

(A) Sampling vectors for {\it homogenized}  coefficient  
$\widetilde{a}_0(x) $, for the oscillatory one $a_\epsilon(x)$, 
as well as for the loading vector ${\bf f}$ all have low QTT ranks.

(B) Numerical implementation of $\mathbb{A}_0^{-1}$ is cheap,  
the solution of ``homogenized'' equation
$ \mathbb{A}_0 {\bf v}_0 = {\bf f}$
has low QTT rank.

(C) For FEM-Galerkin approximation matrix $\mathbb{A}_\epsilon$
the spectral equivalence relation holds
$$
\lambda_0 \mathbb{A}_0 \leq \mathbb{A}_\epsilon \leq \lambda_1 \mathbb{A}_0, 
\quad \lambda_0, \lambda_1 >0.
$$

Assumptions (A) -- (C) are satisfied in all numerical examples presented in the following.

The tensor iterative scheme with QTT rank truncation is described as follows.
Choose the threshold parameter $\delta>0$, 
and denote by ${\cal T}_\delta$ the tensor operation producing almost the best
QTT $\delta$-approximation to the target rank-structured tensor.
Then the exact iteration (\ref{2.17}) on the full finite element space 
is modified as follows:
Starts with $\v_0={\Bbb A}^{-1}_0\f$ presented as the low QTT rank tensor
and then compute $\v_{k+1}$, $k=0,1,2, ...$ 
via fix-point iteration performed in the quantized tensor space $\mathbb{Q}_{2,L}$
and accomplished with $\delta$-rank truncation,
\begin{equation}\label{eqn:Tens_Iter}
 {\bf v}_{k+1} =  
{\cal T}_\delta \left(\beta{\bf v}_{0} - \mathbb{B} {\bf v}_{k}\right),\quad 
\mbox{with} \quad \mathbb{B}= \beta \mathbb{A}_0^{-1} \mathbb{A}_\epsilon - \mathbb{E}.
\end{equation}

Condition (C) ensures the geometric convergence rate for the PCG, Preconditioned
Steepest Descent (PSD) and others accelerated  preconditioned iterations performed 
in the tensor format (\ref{eqn:Tens_Iter}).

Taking into account assumptions (A) -- (C),
we arrive at the efficient preconditioned iterative solver for 
the initial FEM system of equations (\ref{2.16}) discretized on finest grid of size $h$
that resolves all local peculiarities  in the matrix coefficients, i.e. $\epsilon \approx n_0 h$.
The natural choice of the rank truncation parameter might be $\delta=O(h^2)$.

\begin{summary}\label{sum:TensModRed}
 Our model reduction approach introduces $\epsilon$-adapted tensor structured 
 approximation to the 
 initial PDE and to the corresponding preconditioner (in turn, based on certain 
 averaging of the oscillating coefficient) that has low-parametric representation
 as the QTT tensor  and fits well
 the almost periodic structure in the coefficients and in the solution, 
 uniformly in the frequency  parameter $1/\epsilon$. 
 Under certain assumptions on the quality of the quantized 
 tensor approximation  to the input data and the solution, the numerical complexity 
 can be reduced to the logarithmic scale, $O(|\log \epsilon|^q)$. 
 \end{summary}

\section{Numerics: iterative solver with logarithmic complexity}
\label{sec:Numerics}

\subsection{Description of problem classes}\label{ssec:Coef_Classif}

\begin{figure}[htbp]
\centering
\includegraphics[width=5.3cm]{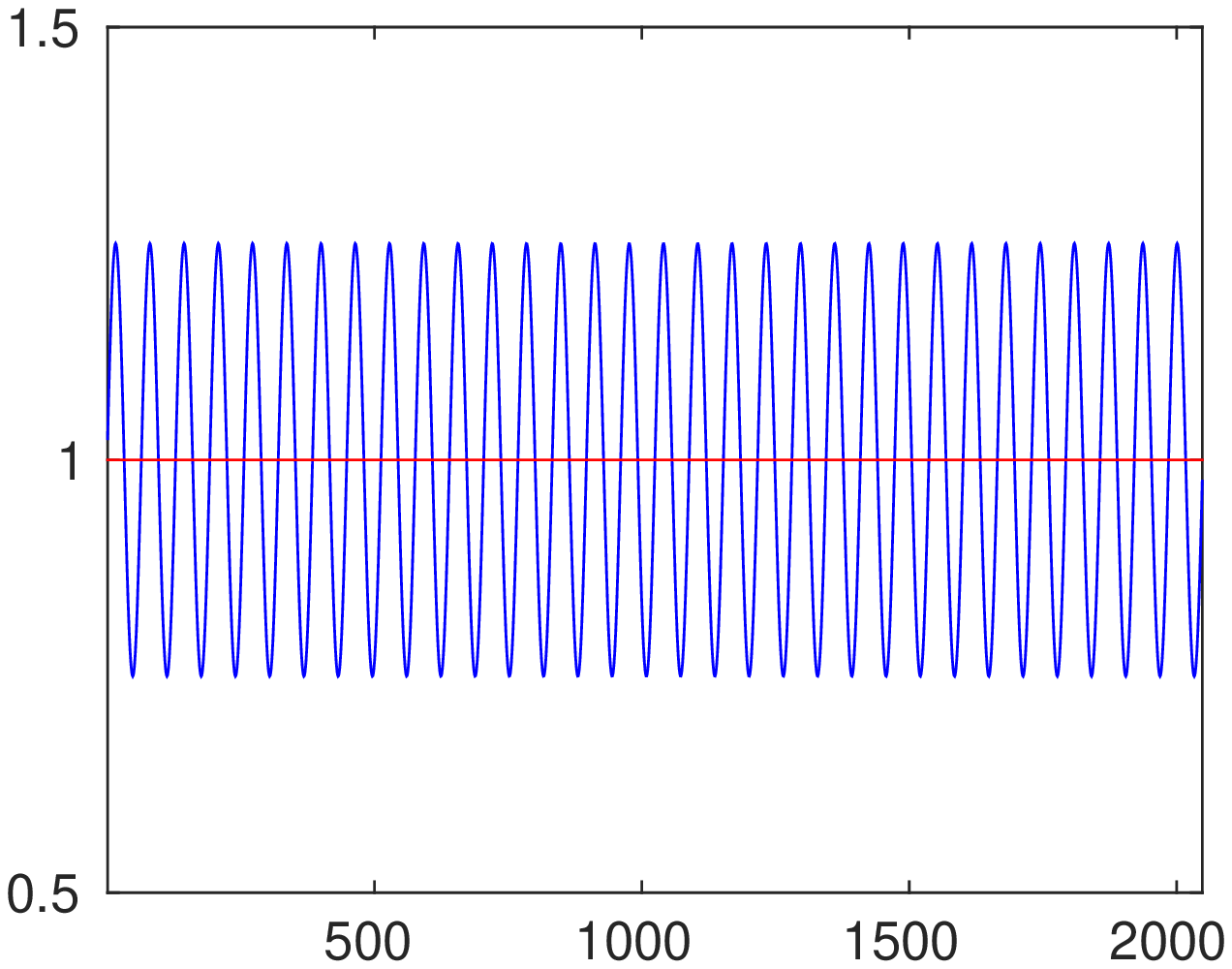}
\includegraphics[width=5.3cm]{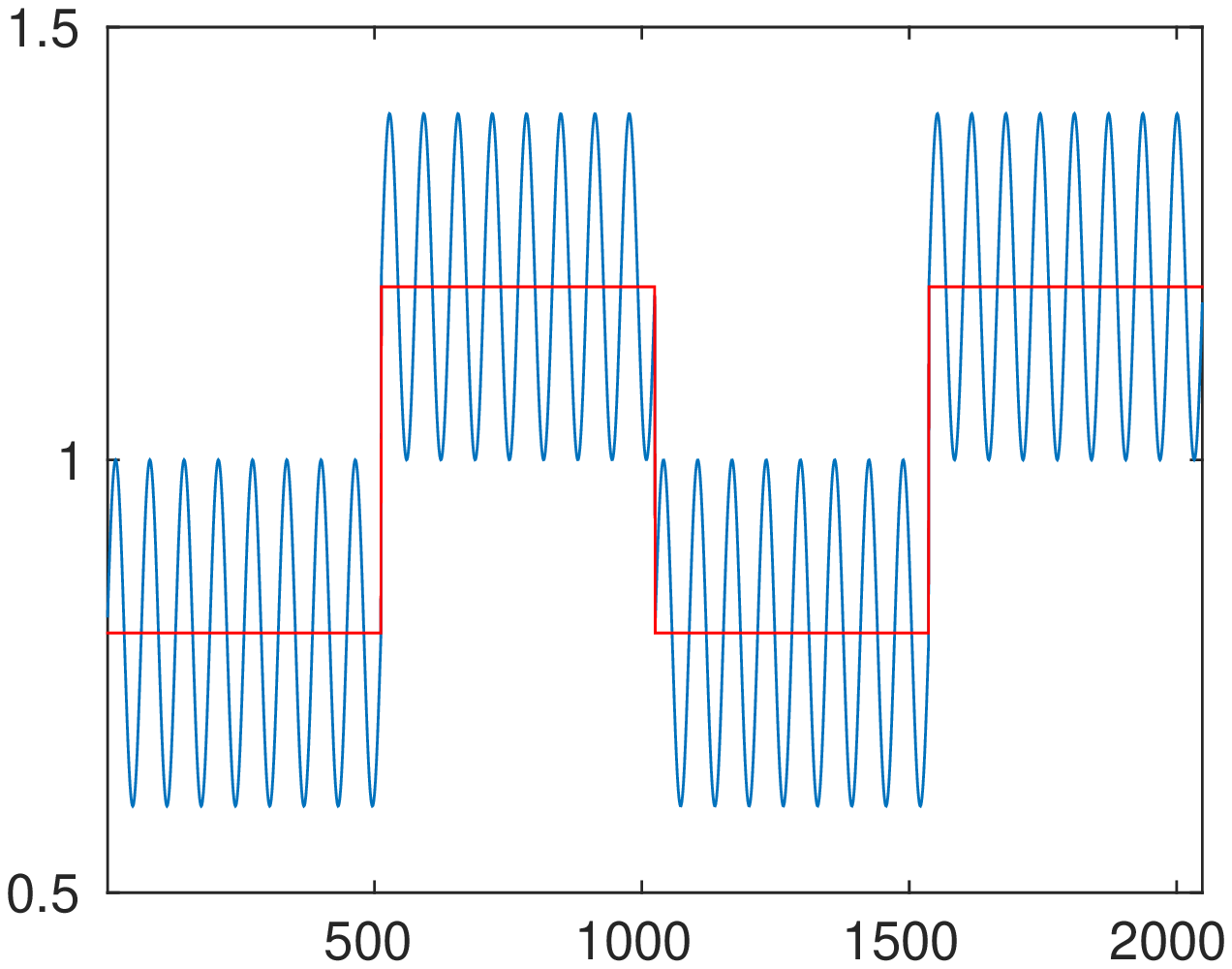}
\includegraphics[width=5.1cm]{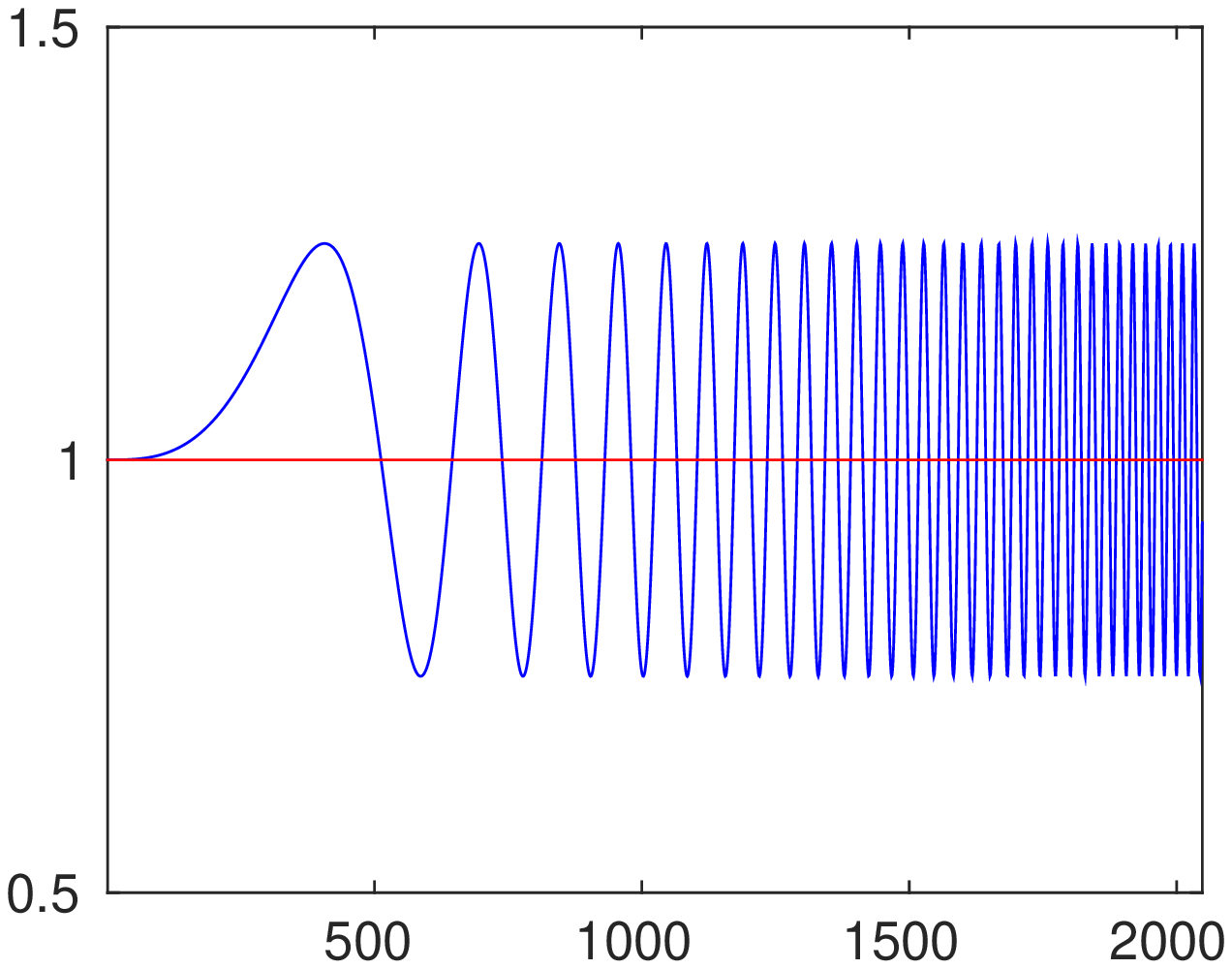}
 \caption{\small Examples of periodic and non-periodic oscillating coefficients.}
\label{fig:Coef1DNperiod}  
\end{figure}

We consider several classes of oscillating or/and jumping coefficients.
Our example of the ideal periodic problem is described by the family of diffusion coefficients
\begin{equation}\label{eqn:PeriodCoef}
 A_\epsilon(x) = C+ \sin (\omega x) > 0, \quad x\in \Omega,
\end{equation}
where the frequency $\omega\in \mathbb{R}$ (i.e. $\epsilon=1/\omega$) may be chosen as
an arbitrarily large constant, see Fig. \ref{fig:Coef1DNperiod}, left.
In this case the rank-$2$ QTT representation of the coefficients vector is suggested
in \cite{KhQuant:09}, see Proposition \ref{prop:KhQuant_Appr:09}, (B). 

The second class of ``modulated periodic'' coefficients is defined by
\begin{equation}\label{eqn:ModPeriodCoef}
 A_\epsilon(x) = C+ g(x)\sin (\omega x),
\end{equation}
where the modulating function $g(x)>0$ should be chosen in such a way that the 
rank of the QTT approximation
to the corresponding function related vector representation of $g(x)$ remains small.
In case of modulated periodic coefficients the QTT rank of the modulated 
function $A_\epsilon(x)$ is bounded by
the product of QTT-ranks for the modulator and oscillator, that is the well known property of
Hadamard product of tensors, see Proposition \ref{prop:KhQuant_Appr:09}, (D).
For the particular choice of the modulator $g(x)$ given by a multi-step function,
see Fig. \ref{fig:Coef1DNperiod}, middle, the QTT-ranks are exactly $2$.
In this case Theorem \ref{thm:QTTrMatr}
ensures that the QTT rank of the resultant equation coefficient ${\bf a}$ does not exceed $4$.

The third class is described by ``exotic'' oscillators which may have nonlinear highly oscillating
behavior and could not be treated by the conventional homogenization methods.
The coefficient is given by
\begin{equation}\label{eqn:exotic_osc}
 A_\epsilon(x) = C+ g(x)\sin (\omega x^m),\quad m=2,3, ...,
\end{equation}
see Fig. \ref{fig:Coef1DNperiod}, right, where $C=1$, $g(x)=1$, and $m=3$.
Clearly, the first two classes of coefficients are the particular cases
of ``exotic'' oscillator in (\ref{eqn:exotic_osc}).

In the general case of ``exotic`` oscillators (\ref{eqn:exotic_osc}) 
the explicit QTT-rank bounds are not known, however, in most examples considered so far
the numerical tests indicate (see e.g. Table \ref{tab:Table_Times}) the low-rank QTT 
approximations with high accuracy (numerical justification).
The rigorous  QTT approximation analysis is possible for some special classes of oscillating
functions, see \cite{KhorVeit:14}.

\subsection{Numerics and comments}\label{ssec:Numer_Comm}

In the following numerical tests  we use the simple preconditioner matrix corresponding to the
constant coefficient $a_0$ defined by the mean value of the initial highly-oscillating 
function $a_\epsilon$, i.e. $a_0=\mean{a_\epsilon}$. 
In this case $\mathbb{A}_0$ is just the scaled 1D Laplacian.

{\small
\begin{table}[htb]
\begin{center}%
\begin{tabular}
[c]{|r|r|r|r|r|r|r|r|}%
\hline
$N=2^L $, iter.   & $2^{13},\, (it)$  & $2^{14},\, (it)$ & $2^{15},\, (it)$ & $2^{16},\, (it)$   & $2^{17},\, (it)$ & $r({\bf a}_\epsilon)$ & $r({\bf u}_\epsilon)$\\
 \hline \hline
$4$-steps coef.    & $3.4, \,(9)$  &$4.3, \,(9)$ & $4.5, \,(9)$ &$6.7, \,(9)$ &$14.3, \,(14)$ & $2.9$ & $4.96$ \\
 \hline
 $C+ \sin (\omega x)$   & $0.97, \,(5)$  & $1.2, \,(5)$  & $1.3, \,(5)$ & $2.0, \,(6)$ & $2.1, \,(6)$ & $2.67$ &  $3.7$ \\
 \hline
$C+ \sin (\omega x^3)$ & $5.3, \,(5)$ & $10.0, \,(6)$  & $9.95, \,(6)$ & $11.98, \,(6)$    & $16.2, \,(5)$ & $7.53$ & $8.24$  \\
 \hline
 \end{tabular}
\caption{CPU times (sec) and rank bounds for the three types of oscillating coefficients. }
\label{tab:Table_Times}
\end{center}
\end{table}
}

We apply the Preconditioned Steepest Descent (PSD) iteration with the QTT-rank truncation 
up to given $\delta>0$ to solve the preconditioned system 
of linear equations represented in the high-dimensional quantized tensor space.
In all the experiments we specify the rank truncation tolerance by $\delta = 10^{-7}$.
The frequency $\omega=1/\epsilon$ was parametrized in the form $\omega=2\pi K$, 
with the particular choice $K=64$. 
The numerical results with smalled or larger frequency parameter $K$ demonstrated the similar
features.

\begin{figure}[htbp]
\centering
\includegraphics[width=5.3cm]{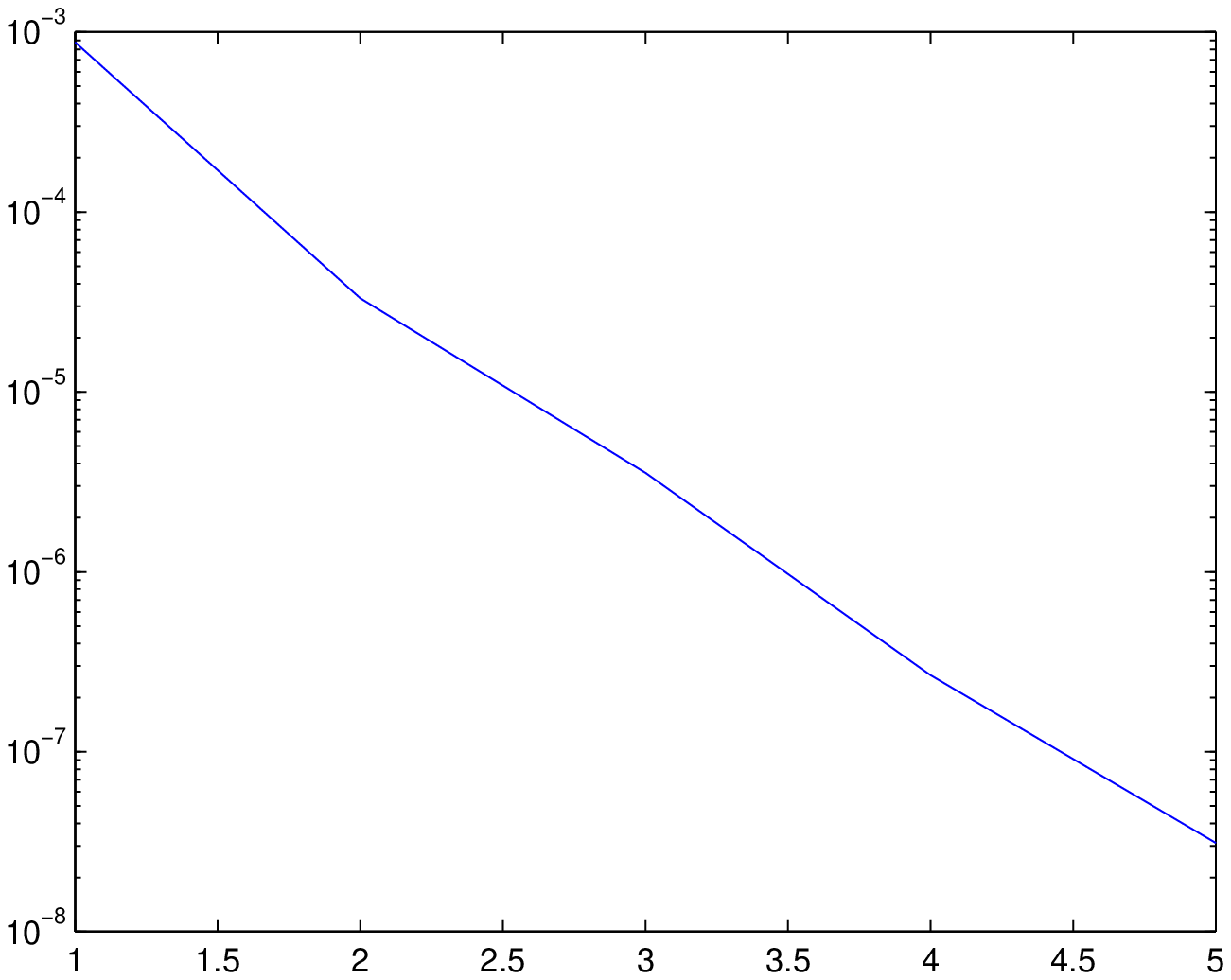}
\includegraphics[width=5.3cm]{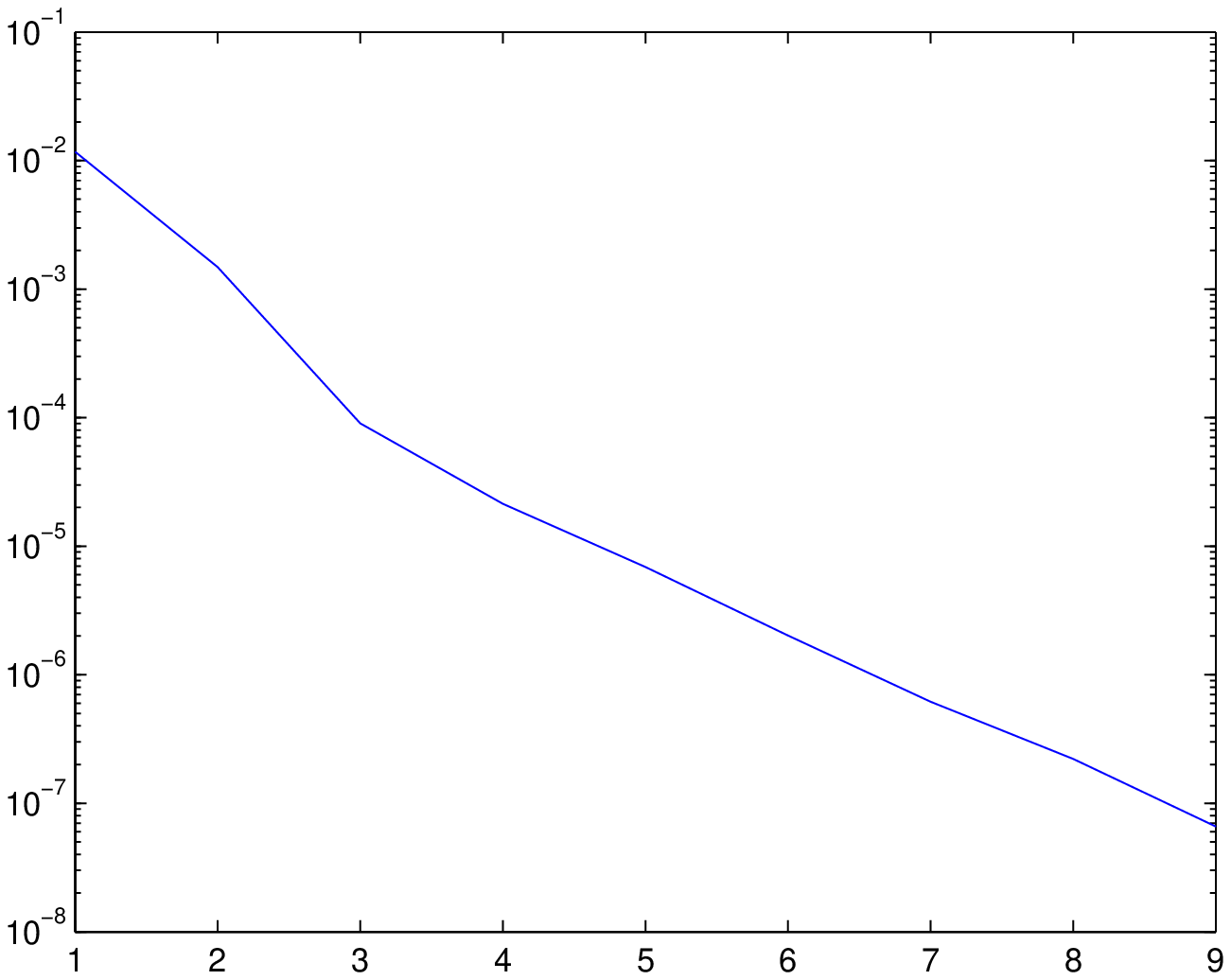}
\includegraphics[width=5.3cm]{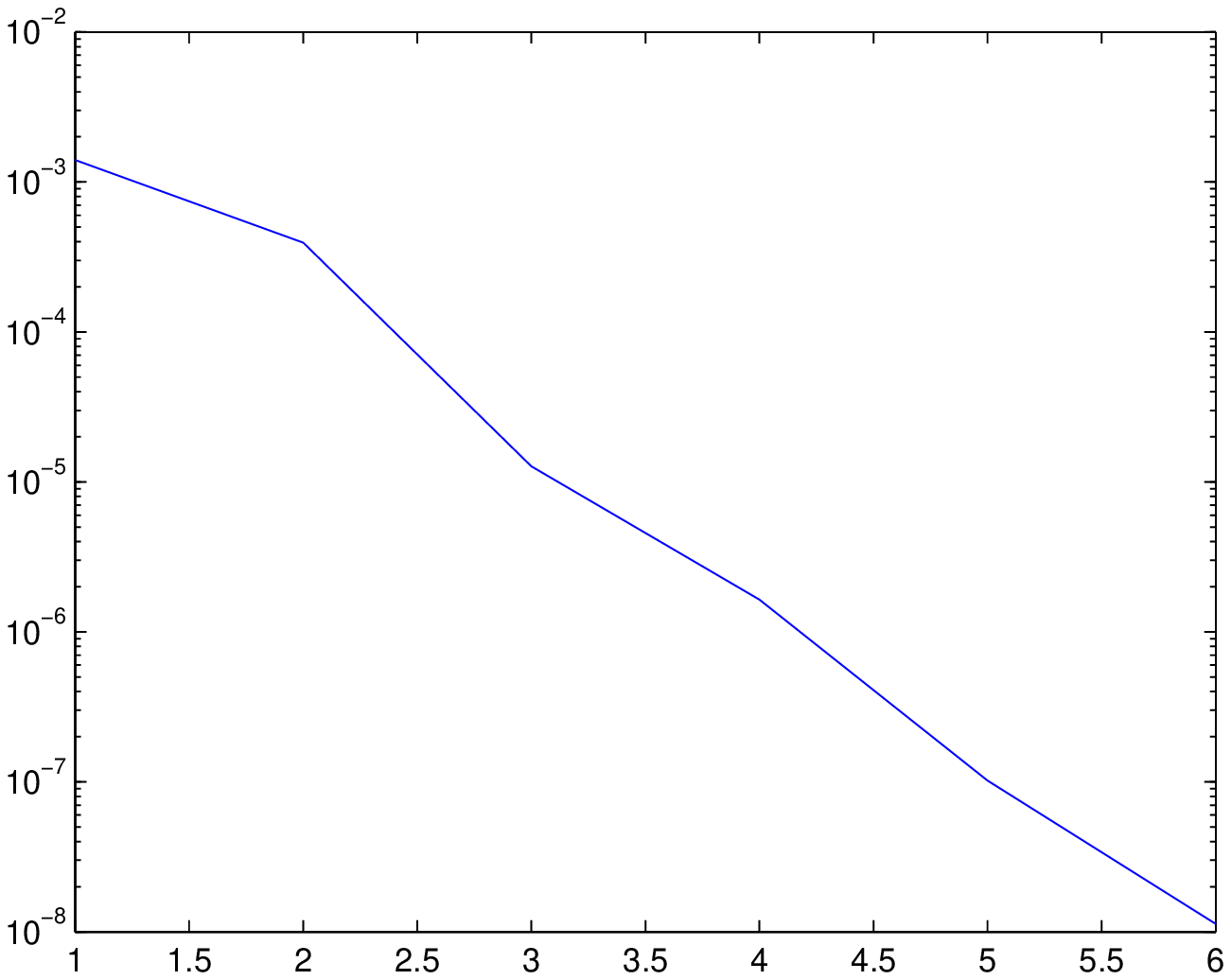}
 \caption{\small PSD iteration history for the periodic (left), jumping oscillating, 
 and cubically oscillating (right) coefficients, see Fig. \ref{fig:Coef1DNperiod}.}
\label{fig:GMRES_conv}
\end{figure}

We observe the uniform geometric convergence of PSD iteration for each of three 
examples considered above, see Table \ref{tab:Table_Times} and Fig. \ref{fig:GMRES_conv}.
Moreover, results in Table \ref{tab:Table_Times} indicate the only logarithmic grows of CPU time 
per iteration with respect to the grid-size $N$. A posteriori error control indicates
the convergence up to $10^{-5}$ in the $H_1$-norm on the finest grid.


The exact solution ${\bf u}_\epsilon$ is calculated with the approximation error of 
the same order $\delta$ (independently of $\epsilon$)
that improves dramatically the approximation error of the homogenized solution.
Notice that 
in all cases presented in Table \ref{tab:Table_Times} the following numerical precision
has been achieved
$\|u_\epsilon - {\bf v}_\delta  \|_0 \backsimeq 10^{-7}$, and
$\|u_\epsilon - {\bf v}_\delta  \|_1 \backsimeq 10^{-6}$.

\begin{figure}[htbp]
\centering
\includegraphics[width=5.4cm]{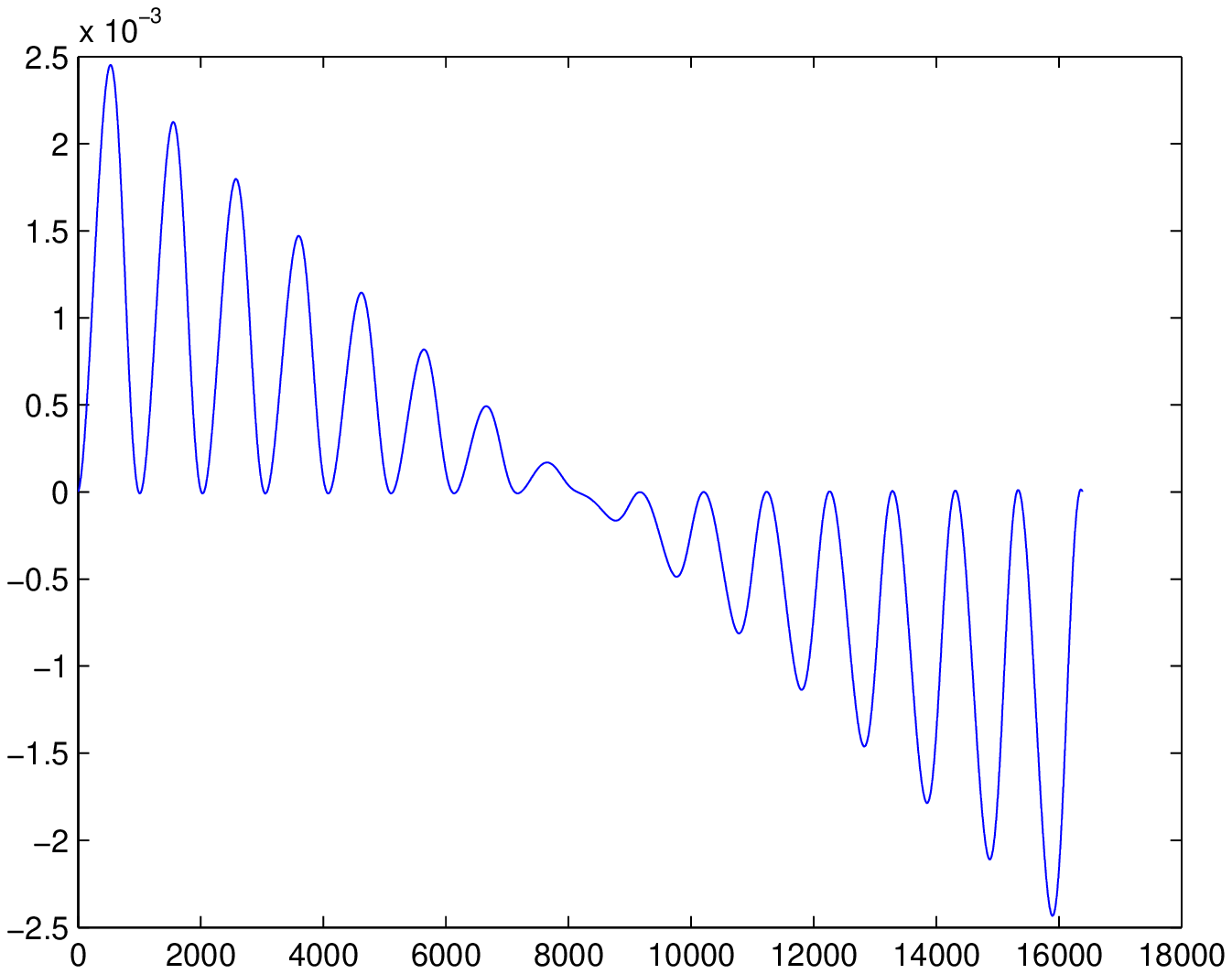}
\includegraphics[width=5.4cm]{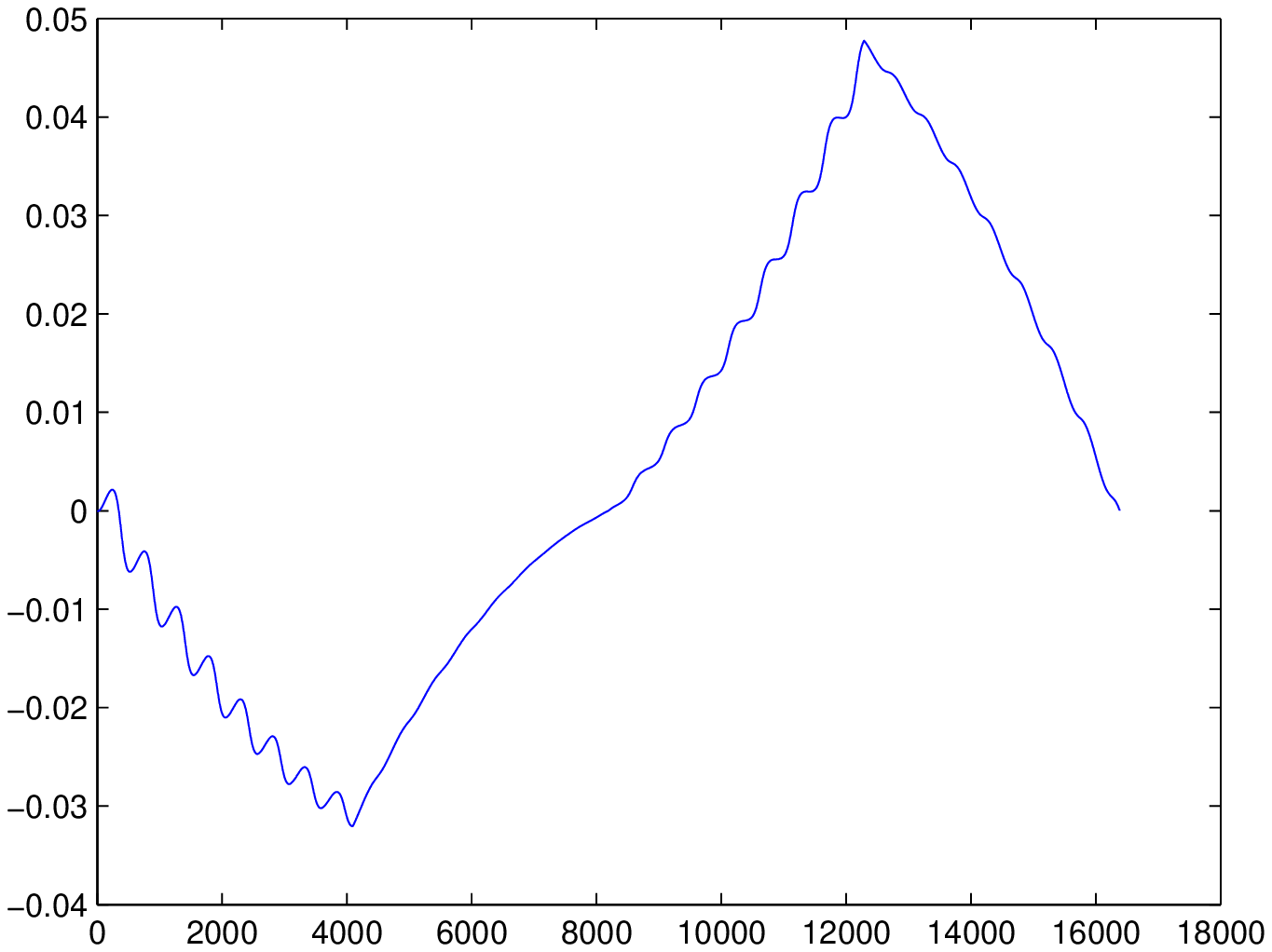}
\includegraphics[width=5.4cm]{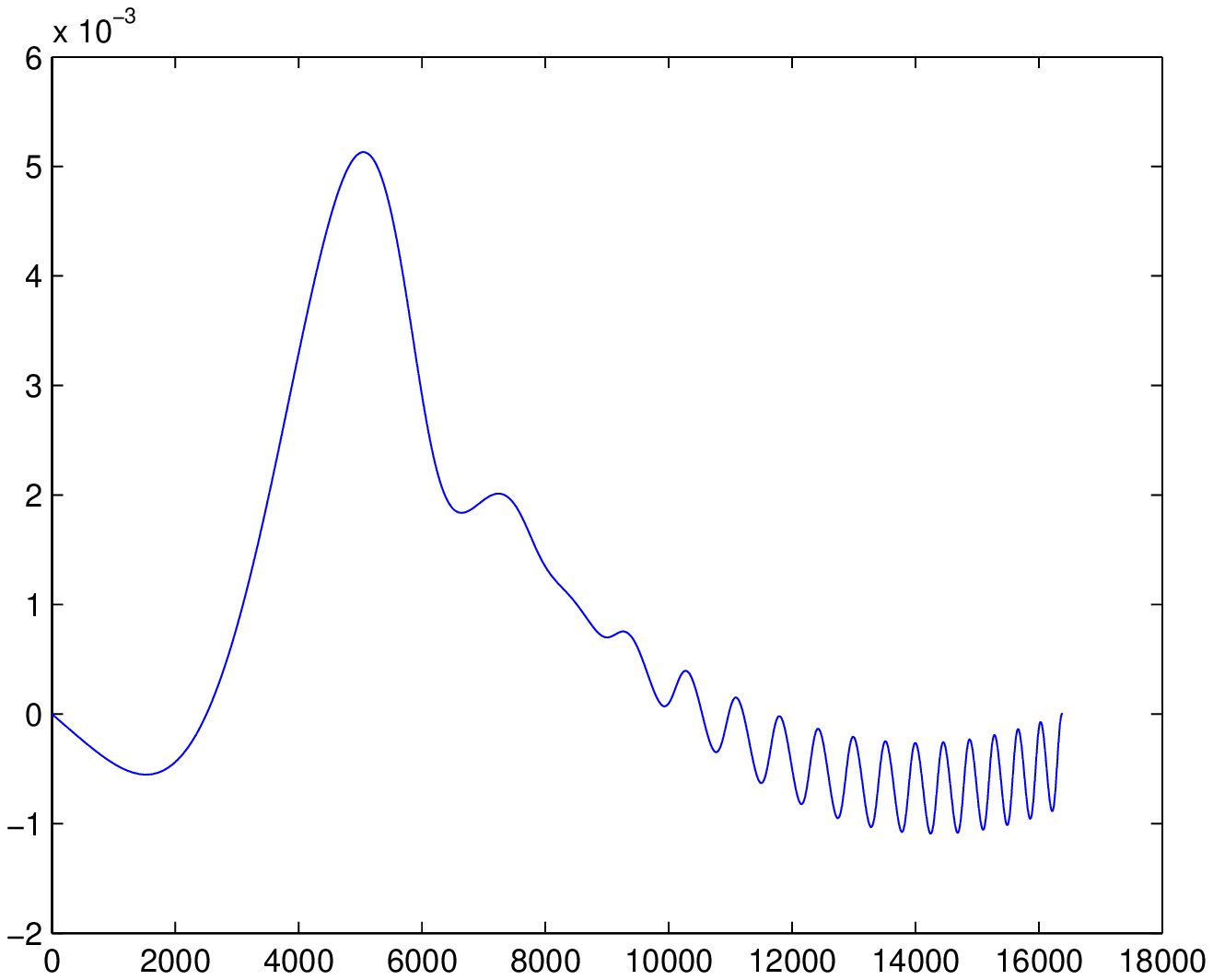}
\caption{\small Difference between exact and homogenized solutions for three 
cases in Fig. \ref{fig:Coef1DNperiod}.}
\label{fig:ErHomExact1DNperiod}
\end{figure}

\begin{figure}[htbp]
\centering
\includegraphics[width=5.4cm]{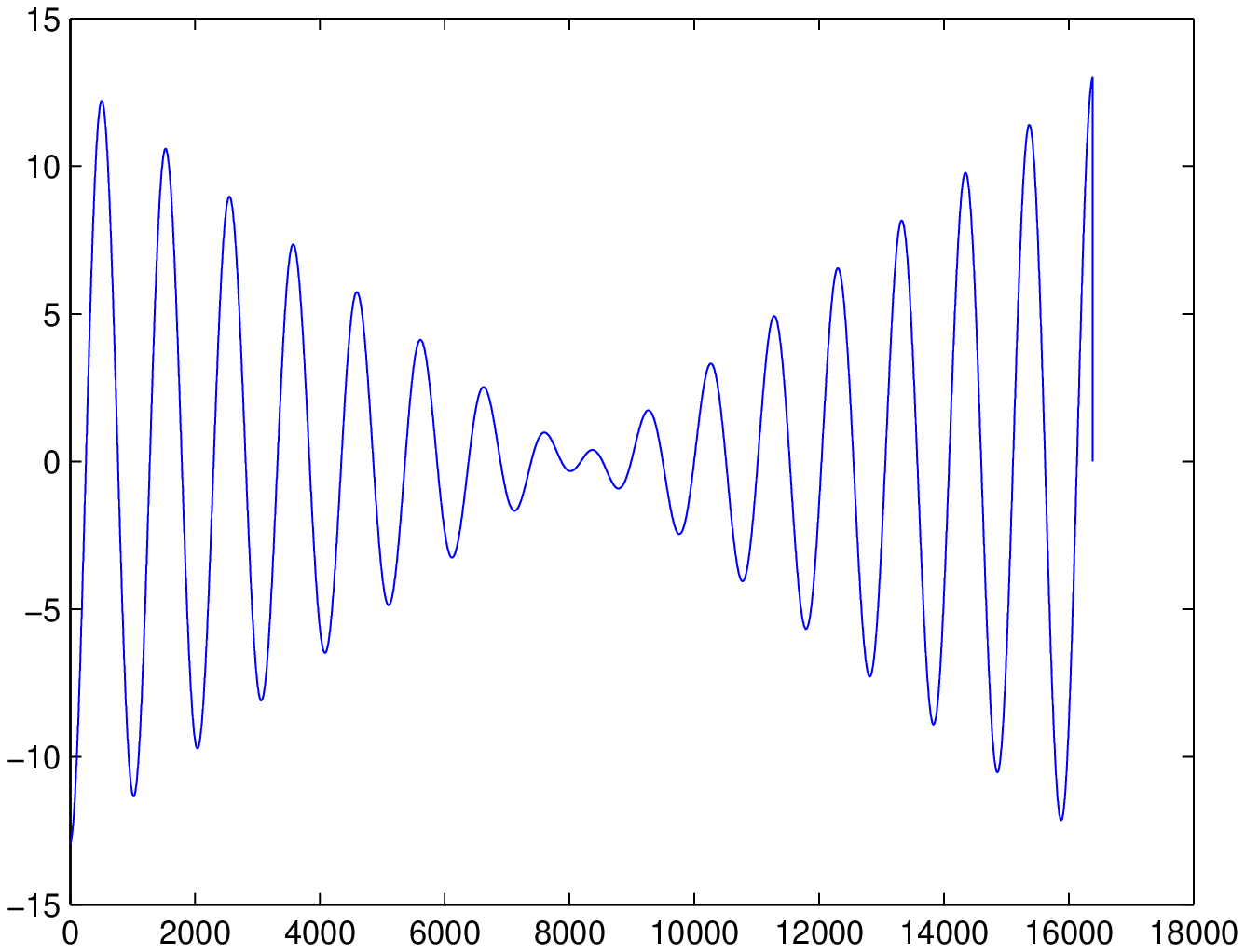}
\includegraphics[width=5.4cm]{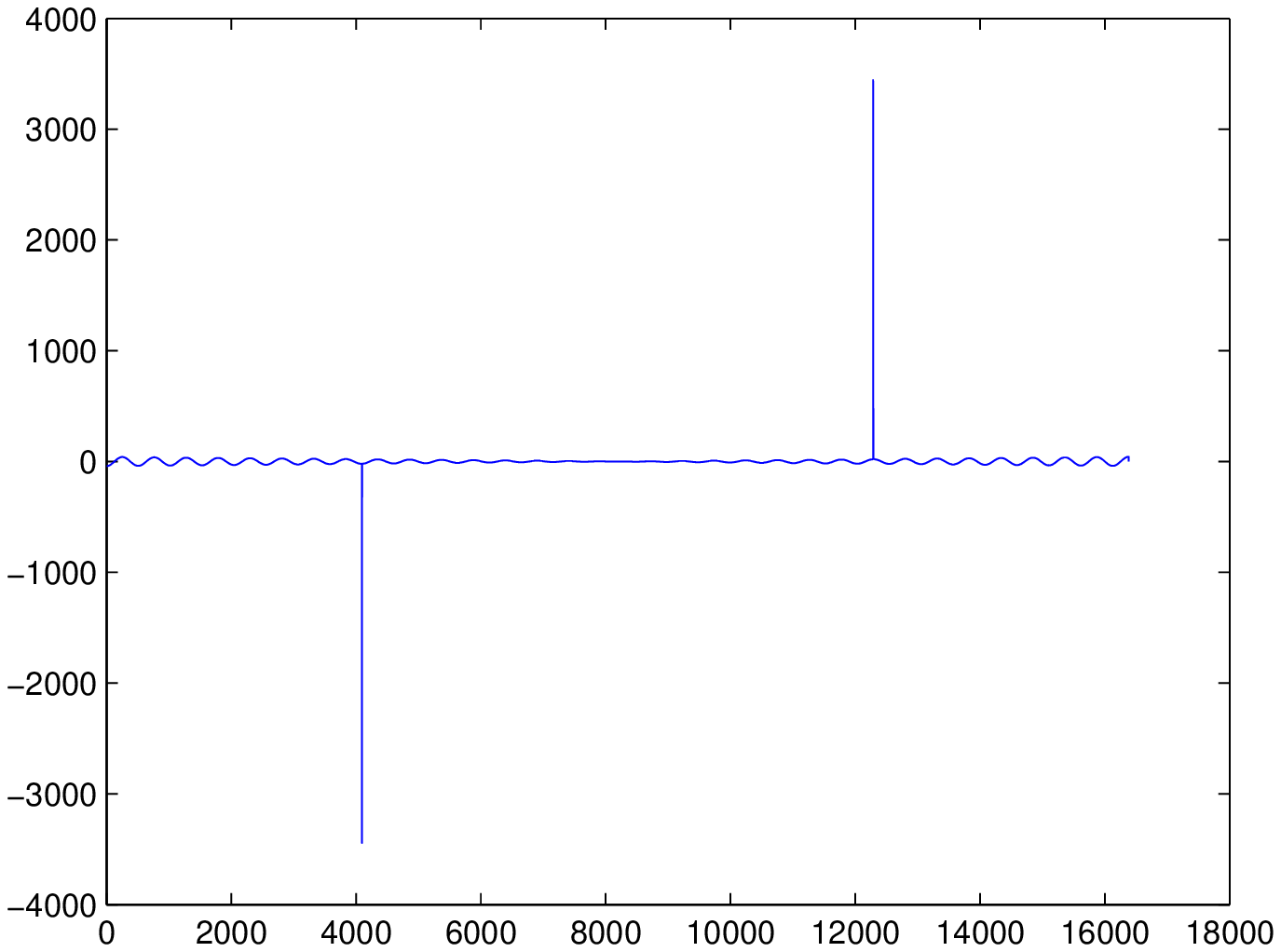}
\includegraphics[width=5.4cm]{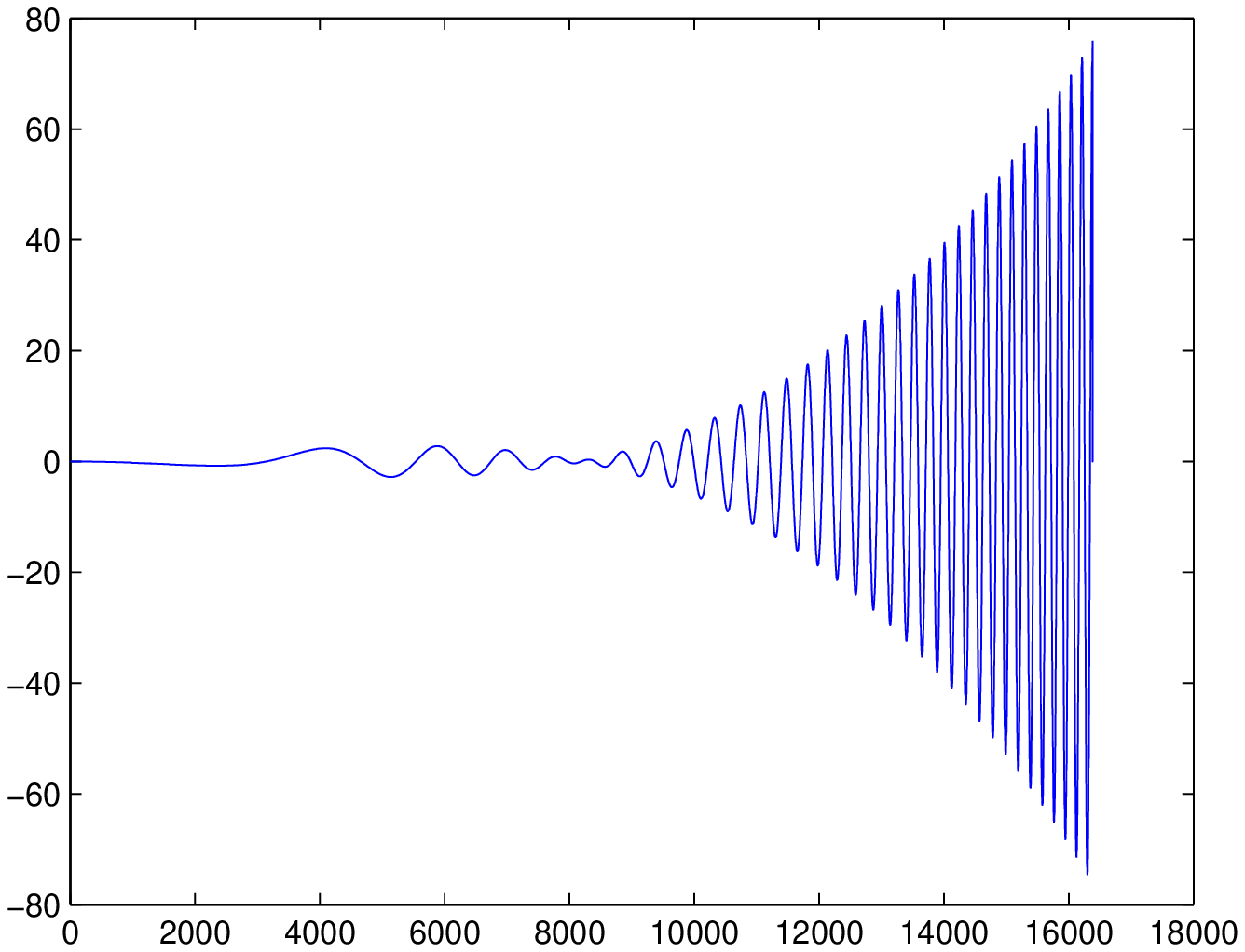}
 \caption{\small Residual for the homogenized periodic (left), $4$-step jumping,
 and cubically oscillating (right) solutions.}
\label{fig:ResHomExact1DNperiod}
\end{figure}

Though the difference between exact and homogenized solutions may be  
of order $O(\sqrt{\epsilon})$, see Fig. \ref{fig:ErHomExact1DNperiod}, 
the residual remains to be large    as indicated in Fig. \ref{fig:ResHomExact1DNperiod},
no convergence in higher derivatives as $\epsilon \to 0$ is observed.

\begin{figure}[htbp]
\centering
\includegraphics[width=6.0cm]{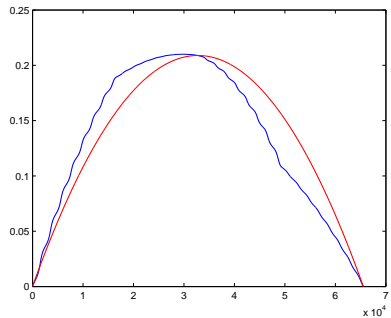}
 \caption{\small Exact (red) and homogenized solutions ($4$-jumping coefficients in
 Fig. \ref{fig:Coef1DNperiod}, middle).}
\label{fig:Hom_vs_Exact1DNperiod}
\end{figure}

In the case modulated periodic coefficients (see Fig. \ref{fig:Coef1DNperiod}, middle, right)
the standard homogenization theory does not provide the convergence even in the limit
of small parameter $\epsilon$. 
In this case Figure \ref{fig:Hom_vs_Exact1DNperiod} demonstrates the systematic error
between the exact (red) and homogenized solutions.

{\bf Acknowledgment.} The second author thanks
Max Planck Institute for Mathematics in the Sciences in Leipzig for support.

{\footnotesize

\end{document}